\documentclass[a4paper]{article}
\usepackage{amssymb}
\usepackage{amsthm}
\usepackage{amsmath}
\usepackage[T1]{fontenc}
\usepackage[utf8]{inputenc}
\usepackage
{hyperref}

\newtheorem{definition}{Definition}[section]
\newtheorem{theorem}{Theorem}[section]
\newtheorem{proposition}{Proposition}[section]

\newtheorem{remark}{Remark}[section]
\newtheorem{corollary}{Corollary}[section]
\newtheorem{example}{Example}[section]

\newcommand\bop {\noindent\textbf{Proof}}
\newcommand\eop {\hbox{}\nobreak\hfill\vrule width 2mm height 2mm depth 0mm
	\par \goodbreak \smallskip}

\newcommand{\integ}[2]{\displaystyle \int_{#1}^{#2}}
\newcommand{\dint}{\displaystyle\int}

\newcommand\R{\mathbb{R}}
\newcommand\Q{\mathbb{Q}}
\newcommand\N{\mathbb{N}}

\newcommand\E{I\!\!E}

\renewcommand\it{\textit}
\renewcommand\ni {\noindent}

\DeclareMathOperator{\essup}{essup}

\title{Doubly Reflected BSDEs With Stochastic Quadratic Growth: Around The Predictable Obstacles}
\author{E. H. Essaky  \qquad M. Hassani \qquad  C. E. Rhazlane\\ \\
	Universit\'{e} Cadi Ayyad\\ Facult\'{e} Poly-disciplinaire\\
	Laboratoire de Modélisation et Combinatoire\\
	D\'{e}partement de Math\'{e}matiques et Informatique\\ B.P. 4162,
	Safi, Maroc.\\  \\ E-mails : essaky@uca.ac.ma\,\, m.hassani@uca.ma \,\, charafeddine.rhazlane@gmail.com}

\date{}

\begin{document}
	
	\maketitle 
	
	\begin{abstract} 
	We prove the existence of maximal (and minimal) solution for one-dimensional generalized doubly reflected backward stochastic differential equation (RBSDE for short) with irregular barriers and stochastic quadratic growth, for which the solution $Y$ has to remain between two \textit{rcll} barriers $L$ and $U$ on $[0,T[$, and its left limit $Y_-$ has to stay respectively above and below two predictable barriers $l$ and $u$ on $]0,T]$. This is done without assuming any $P-$integrability conditions and under weaker assumptions on the input data. In particular, we construct a maximal solution for such a RBSDE when the terminal condition $\xi$ is only ${\cal F}_T-$measurable and the driver $f$ is continuous with general growth with respect to the variable $y$ and stochastic quadratic growth with respect to the variable $z$.
	
	Our result is based on a (generalized) penalization method. This method allow us find an equivalent form to our original RBSDE where its solution has to remain between two new \textit{rcll} reflecting barriers $\overline{Y}$ and $\underline{Y}$ which are, roughly speaking, the limit of the penalizing equations driven by the dominating conditions assumed on the coefficients.

	A standard and equivalent form to our initial RBSDE as well as a characterization of the solution $Y$ as a generalized Snell envelope of some given predictable process $l$ are also given.
	
	\end{abstract}
	
	\ni \textbf{Keys Words:} Doubly reflected backward stochastic differential equation,
Stochastic quadratic growth,
Comparison theorem, Penalization method,
	Snell envelope.
	\medskip
	
	\ni \textbf{AMS Classification}\textit{(2020)}\textbf{: }60H10,
	60H20, 60H05.
	\bigskip

\section{Introduction}

The notion of backward stochastic differential equations with two reflecting barriers
has been first introduced by Civitani\'c and Karatzas \cite{CK}.A
solution for such equation, associated with a coefficient $f$,
a terminal value $\xi$ and two barriers $L$ and $U$, is a quadruple of
processes $(Y, Z, K^+, K^-)$ with values in $\R\times\R^d\times
\R_+\times \R_+$ satisfying:
\begin{equation}
	\label{eq000} \left\{
	\begin{array}{ll}
		(i) &
		Y_{t}=\xi
		+\integ{t}{T}\Big[f(s,Y_{s},Z_{s})ds +dK_{s}^+
		-dK_{s}^- -Z_{s}dB_{s}\Big],\,\, t\leq T,
		\\ (ii)&
		L_t \leq Y_{t}\leq U_{t},\quad \forall t\leq T,\\
		(iii)& \integ{0}{T}( Y_{t}-L_{t}) dK_{t}^+=\integ{0}{T}(
		U_{t}-Y_{t}) dK_{t}^-=0,\,\, \text{a.s.}, \\ (iv) &  K^+, K^-
		\,\,\text{ are continuous non decreasing processes with }\,\, K_0^+
		=K_0^- =0,
	\end{array}
	\right.
\end{equation}
%Due to the cumulative actions of the increasing processes $K^+$ and $K^-$, the solution $Y$ has to remain between the two continuous processes $U$ and $L$.
where  $(B_t)_{t\leq T}$ is the standard Brownian motion. Cvitanic and Karatzas \cite{CK} have proved the existence and uniqueness of solutions if, on the one hand, the coefficient $f$ is uniformly Lipschitz and, on the other hand, the barriers $L$ and $U$ are either regular or satisfy Mokobodski's condition. This later condition essentially postulates the existence of a quasimartingale between the barriers $L$ and $U$. It has also been shown in \cite{CK} that the solution coincides with the value of a stochastic Dynkin game. The link between PDEs with obstacles and RBSDEs has been given in Hamad\`ene and Hassani \cite{HH}. More studies on RBSDEs can be found in \cite{BHM, BY, ekp, GIOQ, HHd, HHO, KLQT, Pexu} and the references therein. %In Hamad\`ene \cite{Ham1}, applications of RBSDEs to
%Dynkin games theory as well as to American game option are given.

The problem of existence of solutions for generalized doubly reflected backward stochastic differential equation (RBSDE for short), which involves an integral with respect to a continuous increasing process, under weaker assumptions on the input data has been studied in \cite{EH1} (see also \cite{PZ,EH2} for the non-reflected case and \cite{Ham} for discontinuous barriers). In \cite{EH1}, the authors have proved the existence of a maximal solution when the terminal condition $\xi$ is ${\cal F}_T-$measurable, the coefficient $f$ is continuous with a general growth with respect to the variable $y$ and a stochastic quadratic growth with respect to the variable $z$. The reflecting barriers $L$ and $U$ are assumed to be continuous. The result has been proved without assuming any $P$-integrability conditions on the terminal data. It should be noted that the integral with respect to the continuous increasing process appears naturally when the authors tried to eliminate the quadratic term by using an exponential transformation. An application of the above result to the Dynkin game problem as well as to the American
game option is studied in \cite{EH3}.

The above result on RBSDE has been generalized by  Essaky, Hassani and Ouknine \cite{EH4} to the case of a RBSDE with two \textit{rcll} barriers which involves a term of the form $\displaystyle{\sum_{t<s\leq
		T}h(s, Y_{s-}, Y_s)}$. %where $h$ is a process with values in $\R$,
% {\color{red}which can be interpreted as a parameter of jump reflection at barriers $L$ and $U$}.
The authors have proved the existence of solutions by establishing first a correspondence between the initial RBSDE and another RBSDE whose coefficients are more tractable. They have shown that the existence of solutions for the initial
RBSDE is equivalent to the existence of solutions for the auxiliary RBSDE. Since the integrability conditions on parameters are weaker, they have made use of approximations and truncations to establish the existence result for the auxiliary RBSDE.

It should be pointed out that, in the case of continuous or \textit{rcll} solutions, all the previous papers on RBSDEs were developed in the framework of continuous or \textit{rcll} obstacles. It is then natural to ask the following question : is there any weaker conditions on the data under which the RBSDE has a solution?

In another context, we recall that, when the barriers is $L^2$-process, Peng and Xu \cite{Pexu} have proved the existence and
uniqueness of solutions for RBSDE under Lipschitz condition on the generator and square integrable data where the following condition has been considered instead of condition $(iii)$ in Equation (\ref{eq000}):
$$
\begin{array}{ll}
	& \forall (L^*, U^*)\in {\cal D}\times {\cal
		D}\quad\text{satisfying}\quad   L_t\leq L_t^* \leq
	Y_{t}\leq U_{t}^* \leq U_{t} \,\,a.e\,\,a.s. \quad \text{we have } \\ &  \integ{0}{T}(
	Y_{t-}-L_{t-}^*)
	dK_{t}^+= \integ{0}{T}( U_{t-}^*-Y_{t-}) dK_{t}^-=0,\,\, \text{a.s.}
\end{array}
$$

%The authors have also shown that this problem is equivalent to find the smallest $g$-supermartingale of BSDE that dominates this obstacle.

%The result has been also generalized to the case of two $L^2-$processes.

In this paper, we are concerned with the study of the  following RBSDE with irregular barriers, for which the solution $Y$ has to remain between two \textit{rcll} barriers $L$ and $U$ on $[0,T[$, and its left limit has to stay respectively above and below two predictable barriers $l$ and $u$ on $]0,T]$, of the form
\begin{equation}
	\label{eq00'} \left\{
	\begin{array}{ll}
		(i) & %\text{ is a f-supersoltion on } [0, T]\text{ with }Y_T = \xi, i.e.
		Y_{t}=\xi
		+\integ{t}{T}\Big[f(s,Y_{s},Z_{s})ds+g(s,
		Y_{s-}, Y_s)dA_s +dK_{s}^+
		-dK_{s}^- -Z_{s}dB_{s}\Big]
		\\ (ii)&
		\forall t\in[0,T[,\,\, L_t \leq Y_{t}\leq U_{t},
		\\ (iii)&  \text{on}\,\,]0,T],\,\, Y_{t-}\leq u_t, \,\,\, d\alpha_t-a.e.\quad l_t\leq Y_{t-},\,\,\, d\delta_t-a.e.\\
		(iv)& \forall (L^*, U^*)\in {\cal D}\times {\cal
			D}\quad\text{satisfying}\quad  \forall t< T, L_t\leq L_t^* \leq
		Y_{t}\leq U_{t}^* \leq U_{t} \quad\text{and}
		\\ &\text{on}\,\,]0,T],\,\,
		U_{t-}^*\leq u_t, \,\,\, d\alpha_t-a.e.\,\,\ \text{and}\,\, l_t\leq
		L_{t-}^*,\,\,\, d\delta_t-a.e.\\ & \text{we have } \integ{0}{T}(
		Y_{t-}-L_{t-}^*)
		dK_{t}^+= \integ{0}{T}( U_{t-}^*-Y_{t-}) dK_{t}^-=0,\,\, \text{a.s.}, \\
		(v)& Y\in {\cal D}, \quad K^+, K^-\in {\cal K}, \quad Z\in {\cal
			L}^{2,d},  \\ (vi)& dK^+\perp  dK^-,
	\end{array}
	\right. \end{equation}

Here, ${\cal D}$ is the set of measurable and right
continuous with left limits (\textit{rcll} for short) processes with values in $\R$, the notation $dK^+\perp  dK^-$ means that $K^+$ and $K^-$ are singular,  $\alpha, \delta$ and $A$ are predictable, right continuous and nondecreasing processes, the generators $f$ and $g$ are assumed to be continuous with general growth with respect to the variable $y$ and $f$ satisfies further the following so-called \textit{stochastic quadratic growth condition}:
$$\forall (y, z)\in
\R\times\R^d,\,\, |f(t, \omega, L_t(\omega)\vee y \wedge
U_t(\omega), z )| \leq \eta_t(\omega)+C_t(\omega)|z|^2,\quad
dtP(d\omega)-\mbox{a.e.}
$$
%The second generator $g$ is assumed to be continuous with general growth with respect to the variable $y$, \it{i.e.}
%$$
%\forall (x, y)\in \R\times \R,\,\,|g(t, \omega, L_{t-}(\omega)\vee x \wedge
%U_{t-}(\omega), L_t(\omega)\vee y \wedge U_t(\omega))| \leq
%\beta_t(\omega),\quad
%A(dt)P(d\omega)-\mbox{a.e.}
%$$

Without assuming any $P-$integrability conditions on the data, we prove the existence of maximal and minimal solutions to Equation (\ref{eq00'}) by constructing two
%using a penalization method. This method allow us to construct
new \textit{rcll} reflecting barriers $\overline{Y}$ and $\underline{Y}$ which are in fact the limit of the penalizing equations driven by the dominating conditions assumed on the data.
With these new \textit{rcll} barriers a new RBSDE is considered.
% equivalent to our initial one
%These barriers will play the role of the dominating \it{rcll} obstacles associated to the solution of another RBSDE equivalent to our initial one.
%However,
Our new idea consists to deduce the solvability of the original RBSDE  (\ref{eq00'}) from the solvability of the new one.
% whose is  ensured by our previous work \cite{EH4}.%, associated with the two \it{rcll} reflecting barriers $\overline{Y}$ and $\underline{Y}$, which is equivalent to our initial RBSDE and its solvability is ensured by \cite{EH4}.

More precisely, we construct two \textit{rcll} reflecting barriers $\overline{Y}$ and $\underline{Y}$ such that: \textit{$(Y, Z, K^+, K^-)$ is a solution of the above RBSDE (\ref{eq00'}) if and only
	if it is a solution of the following RBSDE with the two reflecting \textit{rcll} barriers $\underline{Y}$ and $\overline{Y}$:
	\begin{equation*}
		%\label{eq000'}
		\left\{
		\begin{array}{ll}
			(i) &
			Y_{t}=\xi
			+\integ{t}{T}\Big[f(s,Y_{s},Z_{s})ds+g(s,
			Y_{s-}, Y_s)dA_s+dK_{s}^+ -dK_{s}^--Z_{s}dB_{s} \Big],
			\\ (ii)&
			\forall t\in[0,T[,\,\, \underline{Y}_t \leq Y_{t}\leq
			\overline{Y}_{t},
			\\ (iii)&   \integ{0}{T}(
			Y_{t-}-\underline{Y}_{t-})
			dK_{t}^+=\integ{0}{T} ( \overline{Y}_{t-}-Y_{t-}) dK_{t}^-=0,\,\, \text{a.s.}, \\
			(iv)& Y\in {\cal D}, \quad K^+, K^-\in {\cal K}, \quad Z\in {\cal
				L}^{2,d},  \\ (v)& dK^+\perp  dK^-.
		\end{array}
		\right. \end{equation*}
}
It should be noted that, since the processes $\overline{Y}$ and $\underline{Y}$  are \textit{rcll}, the existence of solutions for this last equation is ensured by the work [\cite{EH4}, Theorem 2.1.].

Our second goal is to write the RBSDE (\ref{eq00'})  in a standard form without introducing the test barriers $L^*$ and $U^*$. To this purpose, we construct two predictable processes $l^{*,\delta}(t)$ and $(-u)^{*,\alpha}(t)$, associated to $l$ and $u$, by the following formulas (see Section \ref{sec2} and the Appendix for more details)
\begin{equation}\label{l}
l^{*,\delta}(t)=	\inf_{n}\bigg\{-nt+ \inf\Big\{ a\in \R : \dint_0^t\Big[l_s(\omega)+ns-
a\Big]^+\,d\delta_s(\omega)=0\Big\}\bigg\},
\end{equation}
\begin{equation}\label{u}
(-u)^{*,\alpha}(t)=	\inf_{n}\bigg\{-nt+ \inf\Big\{ a\in \R : \dint_0^t\Big[-u_s(\omega)+ns-
a\Big]^+\,d\alpha_s(\omega)=0\Big\}\bigg\},
\end{equation}
and then we prove that our RBSDE (\ref{eq00'}) can be written in a standard form as follows
\begin{equation*}
	\left\{
	\begin{array}{ll}
		(i) &  Y_{t}=\xi
		+\integ{t}{T}\left[f(s,Y_{s},Z_{s})ds+g(s,
		Y_{s-}, Y_s)dA_s+dK_{s}^+ -dK_{s}^-
		-Z_{s}dB_{s}\right],
		\\ (ii)&
		\forall t\in]0,T],\,\, L_{t-}\vee l^{*,\delta}(t)\leqslant Y_{t-}\leqslant -(-u)^{*,\alpha}(t)\wedge U_{t-},\\  (iii)&
		\integ{0}{T}( Y_{t-}-[L_{t-}\vee l^{*,\delta}(t)])
		dK_{t}^+= \integ{0}{T}( [U_{t-}\wedge -(-u)^{*,\alpha}(t)]-Y_{t-}) dK_{t}^-=0,\,\, \text{a.s.}, \\
		(iv)& Y\in {\cal D}, \quad K^+, K^-\in {\cal K}, \quad Z\in {\cal
			L}^{2,d},  \\ (v)& dK^+\perp  dK^-.
	\end{array}
	\right. \end{equation*}
This is done by proving the following characterization :
\[\bigg(l_s(\omega)\leqslant Y_{t-}(\omega)\bigg) \;\; d\delta_t(\omega)P(d\omega)\text{-a.e} \quad\text{ if and only if }\quad \bigg(\forall t\in ]0,T],\;\; l^{*,\delta}(t,\omega)\leqslant Y_{t-}(\omega)\bigg)\;\; P\text{-a.s.} \]
A particular form of such RBSDE is related to the notion of generalized Snell envelope of a predictable process $l :=(l_t)_{0\leq t\leq T}$ introduced in \cite{EHO}.
Roughly speaking,  let $l :=(l_t)_{0\leq t\leq T}$ be an ${\cal F}_t$-adapted right
continuous with left limits process with
values in $\R$ of class ${D}[0,T]$, that is the family
$(l_{\nu})_{\nu\in \cal T}$ is uniformly integrable, where $\cal T$
is the set of all ${\cal F}_t$-stopping times $\nu$, such that
$0\leq \nu\leq T$. The Snell envelope ${\cal S}_t(l)$ of $l
:=(l_t)_{0\leq t\leq T}$ is defined as
\begin{equation*}
	\label{snell}
	{\cal S}_{t}\left( l \right) =ess\sup_{\nu\in {\cal T}_{t}}I\!\!E%
	\Big[ l _\nu |{\cal F}_{t}\Big],
\end{equation*}
where ${\cal T}_{t}$ is the set of all stopping times valued between
$t$ and $T$. According to the work of Mertens (see \cite{DM}),
${\cal S}$ is the smallest \textit{rcll}-supermartingale of class ${D
}[0,T]$ which dominates the process $l$, \textit{i.e.},
$$
\forall t\leq T,\,\,  l_t\leq {\cal S}_{t}\left( l \right), \quad P\text{-a.s.}.
$$
Suppose now that the process $l$ is neither of class $D[0; T]$ nor a \textit{rcll} process but just a predictable process. Let $L\in {\cal D}$ and
$\delta\in {\cal K}$ and assume that there exists a local martingale
$M_t = M_0 + \dint_0^t \kappa_s dB_s$ such that $P-$a.s.,
$$ L_t\leq M_t \,\,\mbox{on}\,\, [0, T[ \,\,\mbox{and} \quad l_t\leq
M_{t} \,\,d\delta_t-a.e.\,\,\mbox{on}\,\, [0,
T]\,\,\,\mbox{and}\,\,\, l_T\leq M_T. $$ Then, we prove that the minimal solution $Y$  of the following RBSDE with lower
barriers $L$ and $l$,
\begin{equation*}
	\label{eq001000} \left\{
	\begin{array}{ll}
		(i) &  Y_{t}=L_T+\integ{t}{T}\Big[dK_{s}^+ -Z_{s}dB_{s}\Big],
		\\ (ii)&
		\forall t\in]0,T],\,\, L_{t-}\vee l^{*,\delta}(t)\leqslant Y_{t-},\\  (iii)&
		\integ{0}{T}( Y_{t-}-[L_{t-}\vee l^{*,\delta}(t)])
		dK_{t}^+=0,\,\, \text{a.s.}, \\
		(iv)& Y\in {\cal D}, \quad K^+\in {\cal K}, \quad Z\in {\cal
			L}^{2,d}.
	\end{array}
	\right. \end{equation*}
is the smallest
\textit{rcll} local supermartingale satisfying
$$
\forall t\in [0, T[,\,\,L_t\leq Y_t, \,\,l_t\leq Y_{t-} \,\,d\delta-
a.e., \,\,\,\mbox{on} \,\,\,[0 ,T]\,\,\,\mbox{and}\,\,\, l_T\leq
Y_T.
$$  The process $Y$ (denoted by  $\mathcal{S}_.(L, ld\delta)$) is called the generalized Snell envelope  (see Section \ref{sec} for more details).
% It should be noted that Equation (\ref{0}) can be written in a natural and standard form as follows (see section \ref{sec2} and the appendix for more general details)
%
%\begin{equation}
%\label{00} \left\{
%\begin{array}{ll}
%(i) &  Y_{t}=L_T+\integ{t}{T}\Big[dK_{s}^+ -Z_{s}dB_{s}\Big],
%\\ (ii)&
%\forall t\in]0,T],\,\, L_{t-}\vee l^{*,\delta}(t)\leqslant Y_{t-},\\  (iii)&
%\integ{0}{T}( Y_{t-}-[L_{t-}\vee l^{*,\delta}(t)])
%dK_{t}^+=0,\,\, \text{a.s.}, \\
%(iv)& Y\in {\cal D}, \quad K^+\in {\cal K}, \quad Z\in {\cal
%	L}^{2,d},
%\end{array}
%\right. \end{equation}
As an example, if we assume that $\delta_t =\lambda$ the Lebesgue measure, then $Y=\mathcal{S}_.(L_t1_{\{t<T\}}+\xi1_{\{t=T\}}, ld\lambda)$) the solution of the following RBSDE
$$
\left\{
\begin{array}{ll}
	(i) & %\text{ is a f-supersoltion on } [0, T]\text{ with }Y_T = \xi, i.e.
	Y_{t}=\xi
	+\dint_{t}^{T}\Big[dK_{s}^+ -Z_{s}dB_{s}\Big],\; t\leq T,
	\\ (ii)& l^{*,\lambda}_t \leq Y_{t-}\;\text{ on }\,\, ]0, T],
	\\ (iii)&   \E\dint_{0}^{T}( Y_{t-}-l^{*,\lambda}_t)
	dK_{t}^+=0, \\
	(iv)& Y\in {\cal D}, \quad K^+\in {\cal K}, \quad Z\in {\cal
		L}^{2,d},
\end{array}
\right. $$
is the smallest local super-martingale dominating the predictable process $l$, i.e.
$$
l_t\leq Y_t,\,\,d\lambda- a.e\,\,\text{ and }\,\, \xi\leq Y_T.
$$
A second example in the case of $\delta_t=1_{\{T'\leqslant t\}}$, where $T'$ is a stopping time on $[0, T]$ is also given.

{Let us describe our plan. In Section 2, we introduce the definition of our RBSDE with irregular barriers. Remarks on assumptions and main result of the paper is introduced in Section 3.  In Section 4, we prove the existence of maximal (and minimal) solution of the RBSDE. Section 5 is devoted to find an equivalent and standard form to our initial RBSDE. An application to the notion of generalized Snell envelope is given in Section 6. Finally, most of the material needed in Sections 5 and 6 is given, in a more general setting, in the appendix.}

\section{Definition of a solution for RBSDEs}
\subsection{Notations}
Let $(B_t)_{t\leq T}$ be a  Brownian motion defined on some probability space $(\Omega, {\cal F}, P)$ and let
$({\cal F}_t)_{t\leq T}$ be the standard augmentation of the filtration generated by $(B_t)_{t\leq
	T}$.

% and ${\cal F}_0$ contains all $P$-null sets of $\cal F$.
%
% Note that
%$({\cal F}_t)_{t\leq T}$ satisfies the usual conditions,
%{i.e.} it is right continuous and complete.
For simplicity, we omit sometimes dependence on $\omega$ of some processes or random
functions.
\medskip

%Let us first introduce the following sets:
The following sets will be frequently used in the sequel.

$\bullet$ $\cal P$ is the sigma algebra of $({\cal F}_t)_t$-predictable
sets on $[0,T]\times\Omega .$

$\bullet$ ${\cal D}$ is the set of ${\cal P}$-measurable and right
continuous with left limits (\textit{rcll} for short) processes
$(Y_t)_{t\leq T}$ with values in $\R$.

$\bullet$ For a given process $Y\in {\cal D}$, we denote :
$Y_{t-}=\displaystyle{\lim_{s\nearrow t}Y_s}, t\leq T$
$(Y_{0-}=Y_0)$, and $\Delta_s Y = Y_s - Y_{s-}$ the size of its
jump at time $s$.

$\bullet$ ${\cal K} := \{ K\in  {\cal D}\quad :\quad K \quad\text{is
	nondecreasing and }\,\, K_0  =0\}$.

For a given process $V\in {\cal K} - {\cal K}$ and for
each $\omega\in\Omega$, $dV_t(\omega)$ denotes the signed measure on
$([0, T], {\cal B}_{[0,T]})$ associated to $V_t(\omega)$ where
${\cal B}_{[0,T]}$ is the Borel sigma-algebra on $[0, T]$ and
$$\dint_a^b dV_s = V_b -V_a = \dint_{]a, b]}dV_s.$$

%{\color{red}
%$\bullet$ ${\cal K}^{c} := \{ K\in {\cal K}\quad :\quad \Delta_t K
%=0,\,\, \forall t\in ]0, T]\}$.

%$\bullet$ ${\cal K}-{\cal K}$ is the set of ${\cal P}$-measurable
%and \it{rcll} processes $(V_t)_{t\leq T}$ such that there exist
%$V^+, V^-\in{\cal K}$ satisfying : $V = V^+ -V^-$. In this case, for
%each $\omega\in\Omega$, $dV_t(\omega)$ denotes the signed measure on
%$([0, T], {\cal B}_{[0,T]})$ associated to $V_t(\omega)$ where
%${\cal B}_{[0,T]}$ is the Borel sigma-algebra on $[0, T]$.

%$\bullet$ For a given process $V\in {\cal K}-{\cal K}$, we define :
%$$\dint_a^b dV_s = \dint_{]a, b]}dV_s= V_b -V_a \;\text{ and }\; V_t^c = V_t
%-\displaystyle{\sum_{0 <s\leq t}}\Delta_s V.$$

$\bullet$ ${\cal L}^{2,d}$ the set of $\R^d$-valued and $\cal
P$-measurable processes $(Z_t)_{t\leq T}$ such that
$$\integ{0}{T}|Z_s|^2ds<\infty, P- a.s.$$

\noindent The following notations are also needed :

$\bullet$ For a set $B$, $1_B$ denotes the indicator of $B$.

$\bullet$ For each $(a, b)\in\R^2$, $a\wedge b = \min(a, b)$, $a\vee b = \max(a, b)$, $a^+=a\vee0$ and $a^-=(-a)\vee0$.

$\bullet$ For all $(a, b, c)\in\R^3$ such that $a\leq c$,\,\, $a\vee
b\wedge c =(a\vee
b)\wedge c = a\vee
(b\wedge c)$.

%$\bullet$ Let $K^1$ and $K^2$ be two processes in ${\cal K}$. We note $dK^1 \perp dK^2$ if
%there exists a set $D\in {\cal P}$ such that
%$$
%\E\dint_0^T 1_D(s,\omega) dK^1_s(\omega) = \E\dint_0^T
%(1-1_{D})(s,\omega) dK^2_s(\omega) =0.
%$$
\subsection{The data}

%Before giving the definition of our RBSDE, let us first give the definition of two singular measures.
%\begin{definition}
%Let $K^1$ and $K^2$ be two processes in ${\cal K}$. We say that :
%\begin{enumerate}
%\item $K^1$ and $K^2$ are singular if and only if there exists a set
%$D\in {\cal P}$ such that
%$$
%\E\dint_0^T 1_D(s,\omega) dK^1_s(\omega) = \E\dint_0^T
%1_{D^c}(s,\omega) dK^2_s(\omega) =0.
%$$
% This is denoted by $dK^1 \perp dK^2$.
%\item $dK^1 \leq dK^2$ if and only if for each set $B\in {\cal P}$
%$$
%\E\dint_0^T 1_B(s,\omega) dK^1_s(\omega)\leq \E\dint_0^T
%1_{B}(s,\omega) dK^2_s(\omega), \quad \it{i.e.}\quad K_t^1
%-K_s^1\leq K_t^2-K_s^2,\,\,\, \forall s\leq t\quad P-a.s.
%$$
%In this case $\dfrac{dK^1}{dK^2}$ denotes a ${\cal P}-$measurable
%Radon-Nikodym density of $dK^1$ with respect to $dK^2$ which
%satisfies
%$$
%0\leq \dfrac{dK^1}{dK^2}(s,\omega)\leq 1,\quad dK^2_s
%(\omega)P(d\omega)-a.e. \,\,\text{on}\,\, [0,T]\times \Omega.
%$$
%\end{enumerate}
%\end{definition}

Throughout the paper we need the following data.

\begin{enumerate}
	\item \textit{\textbf{Terminal data :}} $\xi$ is a ${\cal F}_T$-measurable one dimensional random variable.
	\item \textit{\textbf{Lower Barriers :}}
	\begin{enumerate}
		\item[a.]  $l:=\left\{ l_{t},\,0\leq t\leq T\right\}$ is a $ {\cal
			P}-${measurable} process.
		\item[b.]  $L:=\left\{ L_{t},\,0\leq t\leq T\right\}$  is a process
		which belong to ${\cal D}$.
	\end{enumerate}
	\item \textit{\textbf{Upper Barriers :}}
	\begin{enumerate}
		\item[a.]  $u:=\left\{ u_{t},\,0\leq t\leq T\right\}$ is a $ {\cal
			P}-${measurable} processes.
		\item[b.]  $U:=\left\{ U_{t},\,0\leq t\leq T\right\}$ is a process
		which belong to ${\cal D}$.
	\end{enumerate}
	\item \textit{\textbf{Drivers (or generators) : }}
	\begin{enumerate}
		\item[a.]  $f : [0,T] \times \Omega \times \R\times
		\R^d\longrightarrow \R$ is a function  such that for every $(y,z)\in \R\times \R^d$,
		\[(t,\omega)\longmapsto f(t,
		\omega, L_t(\omega)\vee y \wedge U_t(\omega), z)\,\,\text{ is }\,\,
		{\cal P}-\text{measurable}.\]

		\item[b.]  $g :  [ 0,T]\times \Omega \times \R \times \R\longrightarrow \R$ is a
		function such that for any $(x, y)\in \R\times \R$,
		\[(t,\omega)\longmapsto g(t,
		\omega, L_{t-}(\omega)\vee x \wedge U_{t-}(\omega), L_t(\omega)\vee
		y \wedge U_t(\omega))\,\,\text{ is }\,\, {\cal P}-\text{measurable}.\]
	\end{enumerate}
	\item \textit{\textbf{{Processes }:}}  $\alpha$, $\delta$ and $A$ are processes in ${\cal K}$.
\end{enumerate}

\subsection{Definition of a solution}
Before giving the definition of our RBSDE, let us first give the following definition.
\begin{definition}
	Let $K^1$ and $K^2$ be two processes in ${\cal K}$. We say that :
	\begin{enumerate}
		\item $K^1$ and $K^2$ are singular if and only if there exists a set
		$D\in {\cal P}$ such that
		$$
		\E\dint_0^T 1_D(s,\omega) dK^1_s(\omega) = \E\dint_0^T
		1_{D^c}(s,\omega) dK^2_s(\omega) =0.
		$$
		This is denoted by $dK^1 \perp dK^2$.
		\item $dK^1 \leq dK^2$ if and only if for each set $B\in {\cal P}$
		$$
		\E\dint_0^T 1_B(s,\omega) dK^1_s(\omega)\leq \E\dint_0^T
		1_{B}(s,\omega) dK^2_s(\omega), \,\, \it{i.e.}\quad K_t^1
		-K_s^1\leq K_t^2-K_s^2,\,\,\, \forall s\leq t\,\, a.s.
		$$
		%		In this case $\dfrac{dK^1}{dK^2}$ denotes a ${\cal P}-$measurable
		%		Radon-Nikodym density of $dK^1$ with respect to $dK^2$ which
		%		satisfies
		%		$$
		%		0\leq \dfrac{dK^1}{dK^2}(s,\omega)\leq 1,\quad dK^2_s
		%		(\omega)P(d\omega)-a.e. \,\,\mbox{on}\,\, [0,T]\times \Omega.
		%		$$
	\end{enumerate}
\end{definition}
Let us now introduce the definition of our RBSDE for which the solution is constrained to stay between given \textit{rcll} processes $L$ and $U$ and ${\cal P}$-measurable processes $l$ and $u$ (conditions $(ii)$ and $(iii)$). Two nondecreasing processes $K^+$ and $K^-$
are introduced in order to push the solution $Y$ to stay between the barriers in a minimal way.   This
minimality property on $K^+$ and $K^-$ is ensured by the generalized Skorohod conditions (condition $(iv)$) together with the additional constraint  $dK^+\perp  dK^-$  (condition $(vi)$). It should be noted that this orthogonality condition is introduced for the first time in the domain of BSDE in the paper \cite{EH1}.

\begin{definition}\label{def1}
	\begin{enumerate}
		\item A quadruple $(Y,Z,K^+,K^-)\in {\cal D} \times {\cal
			L}^{2,d} \times \in {\cal K} \times {\cal K},$
		is a solution of the RBSDE, associated with the data $(\xi, fdt+ gdA,l\,d\delta+Ldt, u\,d\alpha+Udt)$, if $(i)$, $(ii)$, $(iii)$, $(iv)$ and $(vi)$ of Equation (\ref{eq00'}) are satisfied.

		\item We say that the RBSDE (\ref{eq00'}) has a maximal (resp. minimal)
		solution $(Y_t ,Z_t ,K^+_t , K_t^- )_{t\leq T}$ if for any other
		solution $(Y_t^{'} ,Z_t^{'} ,K^{'+}_t , K_t^{'-} )_{t\leq T}$ we have for all $t \leq T$, $Y_t^{'}\leq Y_t$, $P$-a.s.
		(resp. $Y_t^{'}\geq Y_t$, $P$-a.s.).
	\end{enumerate}
\end{definition}

In order to rewrite equation (\ref{eq00'}) to a more tractable form we  introduce the following set:
$$
	Dom
	=\bigg\{Y\in\mathcal{D}\;:\;
	\E\dint_0^T\left(L_{t-}-Y_{t-}\right)^+dt+\left(l_{t}-Y_{t-}\right)^+d\delta_t+\left(U_{t-}-Y_{t-}\right)^-dt+\left(u_{t}-Y_{t-}\right)^-d\alpha_t=0\bigg\}.
$$

\begin{remark}\label{rem20}
	It should be pointed out that, for every processes $Y, Y'$ and $Y''$ in $\mathcal{D}$,
	%   \begin{enumerate}
	%       \item  if $Y,Y'\in Dom^+$ and for all $t\in [0,T[$, $Y_t\wedge Y'_t\leqslant Y''_t$ then $Y''\in Dom^+$.
	%       \item  if $Y,Y'\in Dom^-$ and for all $t\in [0,T[$, $ Y''_t\leqslant Y_t\vee Y'_t$ then $Y''\in Dom^-$.
	if $Y,Y'\in Dom$ and for all $t\in [0,T[$, $Y_t\wedge Y'_t\leqslant Y''_t\leqslant Y_t\vee Y'_t$ then $Y''\in Dom$.
	%   \end{enumerate}
	%
\end{remark}
\begin{remark}\label{rem21}
	Using the set $Dom$, the RBSDE (\ref{eq00'}) can be written as follows :
	\begin{equation}
		\label{eq0} \left\{
		\begin{array}{ll}
			(i) & Y\in Dom, \quad Z\in {\cal
				L}^{2,d}, \quad K^+, K^-\in {\cal K}\;\hbox{such that}\; dK^+\perp  dK^-,
			\\ (ii)& Y_{t}=\xi
			+\dint_t^T\;\Big[f(s,Y_{s},Z_{s})ds+g(s,
			Y_{s-}, Y_s)dA_s+dK_{s}^+ -dK_{s}^-
			-Z_{s}dB_{s}\Big],
			\\
			(iii)& \forall L^*\in Dom,\quad  \E\dint_0^T\;(
			Y_{t-}-L_{t-}^*)^+\,dK_{t}^++(
		U_{t-}^*-	Y_{t-})\,dK_{t}^-=0.
		\end{array}
		\right. \end{equation}
\end{remark}

\section{Main result}
%\subsection{Assumptions on the data}
We shall need the following assumptions:
\medskip
\begin{enumerate}
	\item \textbf{Assumption $(\textbf{A.1})$  on \boldmath{$f$}:}\\
	There exist two processes $\eta \in L^0(\Omega,
	L^1([0,T], dt, \R_+))$ and $C\in {\cal D}$ such that the driver $f$ satisfies the following conditions:
	\begin{itemize}
		\item[\textbf{(a)}]for all $(y, z)\in
		\R\times\R^d,\,\, |f(t, \omega, L_t(\omega)\vee y \wedge
		U_t(\omega), z )| \leq \eta_t(\omega)+C_t(\omega)|z|^2,\quad
		dtP(d\omega)-$a.e.\\

		\item[\textbf{(b)}] $dtP(d\omega)-$a.e.,\,\, the function $(y, z)\longmapsto f(t,
		\omega,  y, z)$ is continuous.
	\end{itemize}\vspace{0.2cm}
	\item \textbf{Assumption $(\textbf{A.2})$ on \boldmath{$g$}:}\\
	There exists $\beta \in L^0(\Omega, L^1([0,T],
	A(dt), \R_+))$ such that:
	\begin{itemize}
		\item[\textbf{(a)}]
		$A(dt)P(d\omega)-$a.e.,\,\, for every $(x, y)\in \R\times \R,$
		$$
		|g(t, \omega, L_{t-}(\omega)\vee x \wedge
		U_{t-}(\omega), L_t(\omega)\vee y \wedge U_t(\omega))| \leq
		\beta_t(\omega),
		$$
		\item[\textbf{(b)}] $A(dt)P(d\omega)-$a.e. the function
		$$
		(x, y)\longmapsto g(t, \omega,  x , y )\,\,\text{is
			continuous.}
		$$
		\item[\textbf{(c)}] $P-$a.s., for every $(t, x)\in]0, T]\times \R$,\,\, the function
		$$
		y\mapsto y+g(t, \omega, x
		, y )\,\,\Delta_t A\,\,\text{ is
			nondecreasing.}
		$$
	\end{itemize}
	\item \textbf{Assumption $(\textbf{A.3})$:}\\ There exists a semimartingale $S_. = S_0 +
	V^-_. - V^+_. +\dint_0^. \gamma_s dB_s$, with $S_0\in\R, V^{\pm} \in
	\cal K $ and $\gamma\in{\cal L}^{2,d}$, such that $S\in Dom$.
\end{enumerate}

The following theorem constitutes the main result of the paper whose proof is postponed to the next section.
\begin{theorem}\label{thee1}
	If assumptions $(\textbf{A.1})$--$(\textbf{A.3})$ hold then the RBSDE
	(\ref{eq0}) has a maximal and minimal solution.
\end{theorem}
%
%\begin{remark}
%
%\end{remark}
%$\forall t\in
%[0,T[,\,\, L_t \leq S_t\leq U_t$ and \text{on} $]0,T],\,\,
%S_{t-}\leq u_t, \,\,\, d\alpha_t-a.e.\,\,\ \text{and}\,\, l_t\leq
%S_{t-},\,\,\, d\delta_t-a.e.$\\ We assume, without loss of
%generality, that $L_T= U_T = S_T = \xi$.\\
%\subsection{Remarks on the assumptions}\label{rem1}
Let us give the following remarks on the assumptions.

\begin{remark}
\begin{enumerate}
	\item By taking $L_t1_{[0,T[}(t) +\xi1_{\{T\}}(t)$, $U_t1_{[0,T[}(t) +\xi1_{\{T\}}(t)$ and $S_t1_{[0,T[}(t) +\xi1_{\{T\}}(t)$, instead of $L_t$, $U_t$ and $S_t$ respectively, we can assume
	without loss of generality, that
	$$L_T= U_T = S_T = \xi.$$
	\item It should be pointed out that conditions $\textbf{(A.1)(a)}$ and $\textbf{(A.2)(a)}$ hold if the functions $f$ and
	$g$ satisfy the following: $\forall (s,\omega),\,\,\forall x, y\in
	\R,\,\, \forall z\in \R^d$,
	\begin{align*}
		|f(s, \omega, y, z )| &\leq  \varphi(|y|)\;\left(\widetilde{\eta}_s(\omega)\ +\tilde{C}_s(\omega)|z|^2\right),
		\\
		|g(s, \omega,x, y)| &\leq \varphi(|x|+ |y|)\;\widehat{\eta}_s(\omega),
	\end{align*}
	where $\varphi: \R_+
	\longrightarrow \R_+$ is a nondecreasing function, $\widetilde{\eta} \in L^0(\Omega, L^1([0,T], ds,
	\R_+))$, $\tilde{C}\in \mathcal{D}$ and $\widehat{\eta} \in L^0(\Omega, L^1([0,T], dA_s, \R_+))$.

	\noindent Indeed, in conditions $\textbf{(A.1)(a)}$ and $\textbf{(A.2)(a)}$, we have just to take the processes $\eta$, $C$ and $\beta$ as follows:
	\[ \eta_t(\omega)  = \varphi (D_t(\omega))\;\widetilde{\eta}_t(\omega),\;C_t(\omega)  = \varphi (D_t(\omega))\;\tilde{C}_t(\omega)\text{ and }{\beta}_t(\omega)  =\varphi (D_t(\omega))\;\widehat{\eta}_t(\omega),\]
	where $$D_t=2\sup_{s\leqslant t}\Big(U^+_s+L^-_s\Big).$$
	%\begin{align*}
	%\eta_t(\omega) & = \widetilde{\eta}_t(\omega)\varphi (D_t(\omega)),
	%\\  C_t(\omega) & = \tilde{C}_t(\omega)\varphi (D_t(\omega)),
	%\\  {\beta}_t(\omega) & =\widehat{\eta}_t(\omega)\varphi (D_{t-}(\omega)+D_t(\omega)),
	%\end{align*} where $D_t=U^+_t+L^-_t$.
	This means that the functions $f$ and $g$ can have, in particular, a
	general growth with respect to $(x, y)$ and stochastic quadratic growth with respect to $z$. In this respect, assumptions $\textbf{(A.1)(a)}$ and $\textbf{(A.2)(a)}$ are not restrictive.% and one can say more.
	\item Suppose that there exist two processes $\overline{L}, \overline{U}\in Dom$ completely separated, i.e. the processes $\overline{L}$ and $\overline{U}$ are such that:
	$$ \overline{L}_t< \overline{U}_t  \quad \text{on} \quad [0,T[ \quad \text{and}\quad  \overline{L}_{t-}< \overline{U}_{t-} \quad \text{on} \quad ]0,T],
	$$
	then assumption $\textbf{(A.3)}$ holds true. Indeed by setting
	$$
	k_t=1+\sup_{s\leqslant t}\left(|\overline{L}_s|+|\overline{U}_s|\right),\,\, L'_t=\frac{\overline{L}_t}{k_t}1_{(t<T)}+\frac{\overline{L}_{T-}}{k_{T-}}1_{(t=T)}\quad \text{and}\quad U'_t=\frac{\overline{U}_t}{k_t}1_{(t<T)}+\frac{\overline{U}_{T-}}{k_{T-}}1_{(t=T)},
	$$ we have for all $t\in [0,T]$\[-1\leqslant L'_t<U'_t\leqslant1\quad\text{ and }\quad L'_{t-}<U'_{t-}.\] It follows then from [Theorem 4.1., \cite{HHO}] that there is a semimartingale $S'$ such that for all $t\in [0,T]$,
	$$ L'_t\leqslant S'_t\leqslant U'_t,$$ and then the semimartingale $S'_tk_t1_{(t<T)}$ is between $\overline{L}_t1_{(t<T)}$ and $\overline{U}_t1_{(t<T)}$. Hence $S'_tk_t\in Dom$.% (see Remark 2.2 in \cite{EH4}).
\end{enumerate}
\end{remark}

%\item It is not difficult to see that if $L$ or $U$ are
%semimartingales, then $(\bf{A.4})$ holds. Moreover if the barriers
%processes $L$ and $U$ are completely separated on $[0,T[$, \it{i.e.}
%$\forall t\in [0,T[$,\,\,$L_t < U_t$ and $L_{t-} < U_{t-}$ (this is
%equivalent to $\displaystyle{\inf_{0\leq t<T}}(U_t-L_t)>0$), then
%$(\bf{A.4})$ holds. Indeed, let
%$$
%\begin{array}{lll}
%&\beta_t = \displaystyle\sup_{s\leq t} (\mid L_s\mid +\mid U_s\mid )
%\\ & L_t' = \dfrac{L_t}{\beta_t}  1_{\{t<T\}} + \dfrac{L_{T-}}{\beta_{T-}}  1_{\{t=T\}}
%\\ & U_t' = \dfrac{U_t}{\beta_t}  1_{\{t<T\}} + \dfrac{U_{T-}}{\beta_{T-}}
%1_{\{t=T\}}.
%\end{array}
%$$
%Then, $\forall t\in [0,T],\,\, -1\leq L_t'< U_t'\leq 1$ and $
%L'_{t-}< U'_{t-}$. It follows then from the work \cite{HHO} that
%there exists a semimartingale $\overline{S}$ such that $L_t'\leq
%\overline{S}_t\leq U_t',\,\, \forall t\in [0,T]$. Hence, the
%semimartingale $\overline{S_t}\beta_t 1_{\{t<T\}} +\xi 1_{\{t=T\}}$
%is between $L_t$ and $U_t$.

%\subsection{Comparison theorem for maximal solutions }

\section{Proof of the main result}
\noindent This section is devoted to the proof of the existence of maximal solution to Equation (\ref{eq0}) by using a penalization method. This method allows us to construct two \textit{rcll} reflecting barriers which are in fact the limit of penalized equations driven by the dominating conditions assumed on $f$ and $g$. However, our approach consists to deduce the solvability of a RBSDE  (\ref{eq0}) from a suitable RBSDE with two \textit{rcll} reflecting barriers which is equivalent to our initial RBSDE and its solvability is ensured by \cite{EH4}.

\subsection{Comparison theorem for maximal solutions }
Let us now give the following comparison theorem which plays a
crucial rule in the proof of our main result. 
%The proof of this
%comparison theorem is based on an exponential change and an
%approximation scheme, see [Section 3.2. and Lemma 5.1 in \cite{EH4}] for more details. 
For this reason, suppose that assumptions
$\textbf{(A.1)}$--$\textbf{(A.2)}$ hold and $(Y, Z, K^+, K^-)$ is the
maximal solution for the following RBSDE
\begin{equation} \label{Comp}
	\left\{
\begin{array}{ll}
	(i) &  Y_{t}=\xi
	+\integ{t}{T}\left[f(s,Y_{s},Z_{s})ds+g(s,
	Y_{s-}, Y_s)dR_s+dK_{s}^+ -dK_{s}^-
	-Z_{s}dB_{s}\right],
	\\ (ii)&
	\forall t\in[0,T[,\,\, L_t \leq Y_{t}\leq U_{t},\\  (iii)&
	\integ{0}{T}( Y_{t-}-L_{t-})
	dK_{t}^+= \integ{0}{T}( U_{t-}-Y_{t-}) dK_{t}^-=0,\,\, \text{a.s.}, \\
	(iv)& Y\in {\cal D}, \quad K^+, K^-\in {\cal K}, \quad Z\in {\cal
		L}^{2,d},  \\ (v)& dK^+\perp  dK^-.
\end{array}
\right. \end{equation}

Let $(Y', Z', K'^+, K'^-)$ be a solution for
the following RBSDE
\begin{equation} \label{Comp1}\left\{
\begin{array}{ll}
	(i) & Y'_{t}=\xi'
	+\integ{t}{T}\Big[dA'_s +dK'^{+}_{s} -dK'^{-}_{s}
	-Z'_{s}dB_{s}\Big],
	\\ (ii)&
	\forall t\in[0,T[,\,\, L'_t \leq Y'_{t}\leq U'_{t},\\  (iii)&
	\integ{0}{T}( Y'_{t-}-L'_{t-})
	dK'^{+}_{t}= \integ{0}{T}( U'_{t-}-Y'_{t-}) dK'^{-}_{t}=0,\,\, \text{a.s.}, \\
	(iv)& Y'\in {\cal D}, \quad K'^+, K'^-\in {\cal K}, \quad Z'\in
	{\cal L}^{2,d},  \\ (v)& dK'^+\perp  dK'^-,

\end{array}
\right.\end{equation}
  where $R\in {\cal K}, A'\in {\cal K}-{\cal K}$, $L'$ and
$U'$ are two barriers
which belong to ${\cal D}$. \\ \noindent To derive a comparison theorem, we assume the following assumption.\\
\textbf{Assumption $(\textbf{H})$:}
\begin{itemize}
	\item[(H.1)] $\xi' \leq \xi, \,\,$$Y'_t\leq U_t$,\,\, $L'_t\leq Y_t$, $\forall t\in
	[0,T[$.
	\item[(H.2)] $dA'_s\leq f(s, Y'_s, Z'_s)ds+ g(s,Y'_{s-} ,Y'_s)dR_s$ on $]0,T]$.
\end{itemize}
%, $L'$ and
%$U'$ are two barriers
% which belong to ${\cal D}$ such that $L'_t\leq U'_t$, $\forall t\in
% [0,T[$. \\
% \noindent Assume moreover the assumption $\bf{(H)}$:
%\
Suppose that assumptions (A.1)–(A.2) are in force. From Theorem \textbf{6.1} in the paper \cite{EH4}, we have the following comparison theorem for maximal solutions.
\begin{theorem}\label{comparison}[Comparison theorem for maximal solutions]
	\label{th2} Under hypothesis $\textbf{(H)}$, we get
	\begin{enumerate}
		\item $Y'_t\leq Y_t$, for every $t\in [0,T]$, $P-$a.s.
		\item $
		1_{\{U_{t-}=U'_{t-}\}}dK'^{-}_t \leq dK^{-}_t\,\,\,
		\text{and}\,\,\, 1_{\{L_{t-}=L'_{t-}\}}dK^{+}_t \leq dK'^{+}_t. $
	\end{enumerate}
\end{theorem}
\bop. In order to prove Theorem \ref{th2} we should only verify that assumptions of Theorem \textbf{6.1} in \cite{EH4} are satisfied. To begin with, set 
$$
b_t = R_t + |A'|_t.
$$
Let $a \in {\cal K}$, $\alpha \in L^0(\Omega,
L^1([0,T], dt))$ such that 
$$
db_t = \alpha_t\, dt + da_t, \,\,\, da_t\perp dt.
$$
Put 
\begin{equation*}
%	\label{eq00} 
	\left\{
	\begin{array}{ll}
		(i) & \overline{f}(s, y, z) = {f}(s, y, z) + \alpha_s g(s, y, y) \dfrac{dR_s}{db_s} \\ 
		  (ii) &\overline{g}(s, y) = g(s, y, y) \dfrac{dR_s}{db_s}\\ 
	 (iii) 	 & \overline{h}(s, x, y) = g(s, x, y) \dfrac{dR_s}{db_s}\Delta_s a \\ 
	 (iv) 	 & \overline{f'}(s) = \alpha_s \dfrac{dA'_s}{db_s}
	 \\ 
	 (v) 	 & \overline{g'}(s) = \dfrac{dA'_s}{db_s}
	 \\ 
	 (v) 	 & \overline{h'}(s) = \dfrac{dA'_s}{db_s}\Delta_s a.
	\end{array}
\right. \end{equation*}
where $\dfrac{dR}{db}$ (respectively $\dfrac{dA'}{db}$) is the Radon-Nikodym derivative of the measure $dR$ (respectively $dA'$) by the measure $db$.

Hence Equations (\ref{Comp})-(\ref{Comp1}) can be written respectively  as follows : 
\begin{equation} \label{Comp'}
	\left\{
	\begin{array}{ll}
		(i) &  Y_{t}=\xi
		+\integ{t}{T}\left[\overline{f}(s,Y_{s},Z_{s})ds+\overline{g}(s,
	Y_s)da^c_s +dK_{s}^+ -dK_{s}^-
		-Z_{s}dB_{s}\right] + \sum_{t<s\leq T}\overline{h}(s, Y_{s-}, Y_{s}),
		\\ (ii)&
		\forall t\in[0,T[,\,\, L_t \leq Y_{t}\leq U_{t},\\  (iii)&
		\integ{0}{T}( Y_{t-}-L_{t-})
		dK_{t}^+= \integ{0}{T}( U_{t-}-Y_{t-}) dK_{t}^-=0,\,\, \text{a.s.}, \\
		(iv)& Y\in {\cal D}, \quad K^+, K^-\in {\cal K}, \quad Z\in {\cal
			L}^{2,d},  \\ (v)& dK^+\perp  dK^-,
	\end{array}
	\right. \end{equation}
and
\begin{equation} \label{Comp'1}\left\{
	\begin{array}{ll}
		(i) & Y'_{t}=\xi'
		+\integ{t}{T}\Big[\overline{f'}(s) ds +\overline{g'}(s) da_s^c +dK'^{+}_{s} -dK'^{-}_{s} 
		-Z'_{s}dB_{s}\Big] + \sum_{t< s\leq T}\overline{h'}(s),
		\\ (ii)&
		\forall t\in[0,T[,\,\, L'_t \leq Y'_{t}\leq U'_{t},\\  (iii)&
		\integ{0}{T}( Y'_{t-}-L'_{t-})
		dK'^{+}_{t}= \integ{0}{T}( U'_{t-}-Y'_{t-}) dK'^{-}_{t}=0,\,\, \text{a.s.}, \\
		(iv)& Y'\in {\cal D}, \quad K'^+, K'^-\in {\cal K}, \quad Z'\in
		{\cal L}^{2,d},  \\ (v)& dK'^+\perp  dK'^-,
		
	\end{array}
	\right.\end{equation}
where $a^c_s = a_s - \displaystyle{\sum_{0< r\leq s}a_r}$, is the continuous part of the process $a$. By assumptions $(H.1)-(H.2)$ we get 
$$
 \dfrac{dA'_s}{db_s} (\alpha_s ds + da_s)\leq f(s, Y'_s, Z'_s)ds+ g(s,Y'_{s-} ,Y'_s) \dfrac{dR_s}{db_s}(\alpha_s ds + da_s).
$$
Since $ da_s\perp ds$, we obtain from the above inequality that
\begin{equation} \label{Comp''1}\left\{
	\begin{array}{ll}
		(i) &  \dfrac{dA'_s}{db_s} \alpha_s ds \leq \Big[f(s, Y'_s, Z'_s)ds+ \alpha_s g(s,Y'_{s-} ,Y'_s) \dfrac{dR_s}{db_s}\Big]ds \\ 
		(ii) &  \dfrac{dA'_s}{db_s} da_s^c \leq  g(s,Y'_{s-} ,Y'_s) \dfrac{dR_s}{db_s}da_s^c \\  
		(iii) &  \dfrac{dA'_s}{db_s} \Delta_s a \leq  g(s,Y'_{s-} ,Y'_s) \dfrac{dR_s}{db_s}\Delta_s a,\\ 
\end{array}
\right.\end{equation}
which is equivalent to 
\begin{equation} \label{Comp'''1}\left\{
	\begin{array}{ll}
		(i) &  \overline{f'}(s)ds  \leq \overline{f}(s, Y'_s, Z'_s)ds \\ 
		(ii) &   \overline{g'}(s) da_s^c \leq \overline{g}(s, Y'_s)da_s^c   \\  
		(iii) &   \overline{h'}(s) \leq \overline{h}(s, Y'_{s-} ,Y'_s).\\ 
	%	\,\, \forall s\in ]0, T]
	\end{array}
	\right.\end{equation}
It follows that assumptions of Theorem \textbf{6.1} in \cite{EH4} are satisfied for equations (\ref{Comp'}) and (\ref{Comp'1}). Therefore
\begin{enumerate}
	\item $Y'_t\leq Y_t$, for every $t\in [0,T]$, $P-$a.s.
	\item $
	1_{\{U_{t-}=U'_{t-}\}}dK'^{-}_t \leq dK^{-}_t\,\,\,
	\text{and}\,\,\, 1_{\{L_{t-}=L'_{t-}\}}dK^{+}_t \leq dK'^{+}_t. $
\end{enumerate} \eop

%This allow us to apply Theorem \textbf{6.1} in \cite{EH4} to Equations (\ref{Comp'}) and (\ref{Comp'1}) to get the result since, from (\ref{Comp'''1}), assumptions of this theorem are satisfied. Theorem \ref{comparison} is then proved. \eop 
\subsection{Penalization method}
To begin, we set for every
$s\in [0, T]$,
$$
m_s =1+8\, \displaystyle{\sup_{r\leq s}|C_r|}.$$
We consider
$(\underline{Y}^n,\underline{Z}^n,\underline{K}^{n+},\underline{K}^{n-})$
the minimal solution of the following penalized RBSDE where the driver is derived from the conditions assumed on $f$ and $g$ and the semimartingale $S$:
\begin{equation}
	\label{eq00} \left\{
	\begin{array}{ll}
		(i) &  \underline{Y}^n_{t}=\xi
		-\integ{t}{T}\Big[\big(\eta_s+ 4C_s |\gamma_s|^2 +\frac12
		m_s|\underline{Z}^n_s -\gamma_s|^2)ds
		-\big(dV^+_s+dV^-_s+\beta_s dA_s\big)
		\\
		&\qquad\qquad\qquad\quad+n(l_s
		-\underline{Y}^n_{s-})^+d\delta_s+d\underline{K}_{s}^{n+}
		-d\underline{K}_{s}^{n-}
		-\underline{Z}_{s}^ndB_{s}\Big]\,, t\leq T,
		\\ (ii)&
		\forall t\in[0,T[,\,\, L_t \leq \underline{Y}^n_{t}\leq S_{t},\\
		(iii)& \integ{0}{T}( \underline{Y}^n_{t-}-L_{t-})
		d\underline{K}_{t}^{n+}= \integ{0}{T}( S_{t-}-\underline{Y}^n_{t-}) d\underline{K}_{t}^{n-}=0,\,\, P-\text{a.s.}, \\
		(iv)& \underline{Y}^n\in {\cal D}, \quad \underline{K}^{n+},
		\underline{K}^{n-}\in {\cal K}, \quad \underline{Z}^n\in {\cal
			L}^{2,d},
		\\ (v)& d\underline{K}^{n+}\perp  d\underline{K}^{n-},
	\end{array}
	\right. \end{equation}
where $S, V^+, V^-$ and $\gamma$ are the processes appeared
in Assumption $(\textbf{A.3})$.

Let also $(\overline{Y}^n,\overline{Z}^n,
\overline{K}^{n+}, \overline{K}^{n-})$ be the maximal solution of
the following penalized RBSDE:
\begin{equation}
	\label{eq1} \left\{
	\begin{array}{ll}
		(i) &  \overline{Y}^n_{t}=\xi
		+\integ{t}{T}\Big[\big(\eta_s+ 4C_s |\gamma_s|^2 +\frac12
		m_s|\overline{Z}^n_s -\gamma_s|^2\big)ds
		+\big(dV^+_s+dV^-_s+\beta_s dA_s\big)
		\\
		&\qquad\qquad\qquad\quad
		-n(\overline{Y}^n_{s-} -u_s)^+d\alpha_s+d\overline{K}_{s}^{n+}
		-d\overline{K}_{s}^{n-}
		-\overline{Z}_{s}^ndB_{s}\Big]\,, t\leq T,
		\\ (ii)&
		\forall t\in[0,T[,\,\, S_t \leq \overline{Y}^n_{t}\leq U_{t},\\
		(iii)& \integ{0}{T}( \overline{Y}^n_{t-}-S_{t-})
		d\overline{K}_{t}^{n+}= \integ{0}{T}( U_{t-}-\overline{Y}^n_{t-}) d\overline{K}_{t}^{n-}=0,\,\, P-\text{a.s.}, \\
		(iv)& \overline{Y}^n\in {\cal D}, \quad \overline{K}^{n+},
		\overline{K}^{n-}\in {\cal K}, \quad \overline{Z}^n\in {\cal
			L}^{2,d},  \\ (v)& d\overline{K}^{n+}\perp  d\overline{K}^{n-},
	\end{array}
	\right. \end{equation}  We should point out here that, since the barriers are \textit{rcll} and the drivers are of stochastic quadratic growth, the existence of minimal
(resp. maximal) solution to (\ref{eq00}) (resp. (\ref{eq1})) is ensured by the work [\cite{EH4}, Theorem 2.1.].

\subsection{Study of the penalized equations (\protect\ref{eq00})-(\protect\ref{eq1})}
In this subsection, we will prove that limiting processes $\underline{Y}$ and $\overline{Y}$ of $\underline{Y}^n$ and $\overline{Y}^n$ respectively are in $Dom$. %\textit{rcll} satisfying moreover
%$$
% \overline{Y}_{t-}\leq u_t,\,\,\, d\alpha_t-a.e.,\quad \text{on}\quad
%]0,T] \quad \text{and}\quad
%l_t\leq \underline{Y}_{t-}, \,\,\, d\delta_t-a.e.,\quad \text{on} \quad ]0,T].
%$$
To begin with, we recall that the semimartingale $S$ is given by
\[S_t = \xi  - \dint_t^T dV^-_t +\dint_t^T dV^+_t -\dint_t^T \gamma_s dB_s.\]
and consider the solution $\overline{Y}^n$ of Equation \ref{eq1}. Then  assumptions $(H.1)$ and $(H.2)$ of Theorem 4.1 are satisfied by taking 
\begin{align*}
    & \xi'=\xi, L' =S, U'=U,  A' =  V^+- V^-\\
    &  Y'= S, Z'=\gamma, dK^{'+}=dK^{'-} =0,
\end{align*}  and 
\begin{align*}
    & f(s, y, z) = \eta_s+ 4C_s |\gamma_s|^2 +\frac12
    m_s|{z} -\gamma_s|^2, \quad  dR_s = dV^+_s+dV^-_s+\beta_s dA_s+d\alpha_s, \\& g(s,x,y)=-n(x-u_s)^+\frac{d\alpha_s}{dR_s}+\frac{dV^+_s+dV^-_s+\beta_s dA_s}{d\alpha_s} ,
    \quad  L=S, \\
    & Y =  \overline{Y}^n ,  K^{\pm}=\overline{K}^{n\pm}, Z= \overline{Z}^n.
\end{align*}

Applying comparison theorem (Theorem 4.1) to
$Y=\overline{{Y}}^{n}$ and $Y'=S$, it follows that 
$$d\overline{K}_s^{n+}=\mathbf{1}_{\left\{L_{s-}=L'_{s-}\right\}}d\overline{K}_s^{n+}
\leq d{K}_s^{'+} =0.$$
Henceforth
$$d\overline{K}^{n+}= 0.$$

\noindent By a symmetric argument, it follows also that for
every
$n\in\N$,
$$d\underline{K}^{n-}= 0.$$
Again, by using comparison
theorem (Theorem \ref{th2}) we get also that
\begin{equation}\label{equ10'}
	L_t \leq \underline{Y}^n_{t}\leq \underline{Y}^{n+1}_{t}\leq S_{t}
	\leq \overline{{Y}}^{n+1}_{t}\leq \overline{{Y}}^{n}_{t}\leq U_t.
\end{equation}
Set
\begin{equation}\label{equ10}
	\begin{array}{ll}&\overline{Y}_t = \displaystyle{\inf_{n}}\overline{{Y}}^{n}_{t}$,
		\text{ and } $\overline{Y}_t^{\,-} =
		\displaystyle{\inf_{n}}\overline{{Y}}^{n}_{t-},
		\\ &
		\underline{Y}_{t}=\displaystyle{\sup_{n}}\underline{Y}^n_{t}
		\text{ and }
		\underline{Y}_{t}^{\,-}=\displaystyle{\sup_{n}}\underline{Y}^n_{t-}.
	\end{array}
\end{equation}
By letting $n$ to infinity in (\ref{equ10'}) and using assumption
$(\textbf{A.3})$ we get that the semimartingale $S$ is between $\underline{Y}$ and  $\overline{{Y}}$. More precisely,  we have the following.
\begin{proposition}\label{rem0}
	%\begin{enumerate}
	%\item
	For every $t\in [0, T]$, we get
	$$
	L_t \leq \underline{Y}_{t}\leq S_{t} \leq \overline{{Y}}_{t}\leq
	U_t\quad \text{and}\quad L_{t-} \leq \underline{Y}_{t}^{-}\leq
	S_{t-} \leq \overline{{Y}}^{\,-}_{t}\leq U_{t-}.
	$$
	%\item \text{On} $]0,T]$, we have
	%$$
	%l_t\leq \overline{{Y}}_{t}^{\,\,-}, \,\,\, d\delta_t-a.e.\,\,\
	%\text{and}\,\, \underline{Y}_{t}^{-}\leq u_t,\,\,\, d\alpha_t-a.e.
	%$$
	%\end{enumerate}
\end{proposition}
\begin{proposition}\label{pro1} The processes $\overline{Y}$ and $\underline{Y}$ defined by (\ref{equ10}) are in $Dom$. In particular, $\overline{Y}$ and \underline{Y} are \textit{rcll}.

	%, \it{i.e.} they satisfy the
	%following properties:
	%\begin{enumerate}
	%\item $\overline{Y}$ and \underline{Y} are \it{rcll}.
	%\item $\overline{Y}_{t-}\leq u_t,\,\,\, d\alpha_t-a.e.,$ \text{on}
	%$]0,T]$.
	%\item $l_t\leq \underline{Y}_{t-}, \,\,\, d\delta_t-a.e.,$ \text{on}
	%$]0,T]$.
	%\end{enumerate}
\end{proposition}
\bop. Let $R_t= \integ{0}{t}\bigg[\Big(\eta_s+ 4C_s
|\gamma_s|^2\Big)ds + 2dV^-_s+\beta_s dA_s\bigg]$. We have
$$
\begin{array}{ll}
	& \overline{Y}_{t}^n -S_t\\ & = \dint_t^T dR_s + \dint_t^T
	\frac{m_s}{2}\Big|\overline{Z}_{s}^n -\gamma_s\Big|^2ds - n\dint_t^T
	\Big(\overline{Y}_{s-}^n -u_s\Big)^+d\alpha_s - \dint_t^T
	d\overline{K}_{s}^{n-}- \dint_t^T\Big(\overline{Z}_{s}^n -\gamma_s\Big)dB_s.
\end{array}
$$
Then
$$
\begin{array}{ll}
	& m_t\Big(\overline{Y}_{t}^n -S_t\Big) \\& = \dint_t^T m_{s-} dR_s +
	\frac12\dint_t^T \Big|m_s\Big(\overline{Z}_{s}^n -\gamma_s\Big)\Big|^2ds -
	n\dint_t^T m_{s-} \Big(\overline{Y}_{s-}^n -u_s\Big)^+d\alpha_s - \dint_t^T
	m_{s-} d\overline{K}_{s}^{n-} \\ &\qquad\qquad\qquad\quad -
	\dint_t^T m_{s}\Big(\overline{Z}_{s}^n -\gamma_s\Big)dB_s-
	\dint_t^T\Big(\overline{Y}_{s}^n -S_s\Big)dm_s.
\end{array}
$$
Let $\psi(r) = e^{r}$. By  It\^{o}'s formula we have

\[\begin{array}{ll}
	&\psi\Big(m_t(\overline{Y}_{t}^n -S_t)\Big) \\&= 1+\dint_t^T
	\psi'\Big(m_{s-}(\overline{Y}_{s-}^n -S_{s-})\Big) m_{s-} dR_s +
	\frac12\dint_t^T \psi'\Big(m_{s}(\overline{Y}_{s}^n -S_{s})\Big)
	\Big|m_s\Big(\overline{Z}_{s}^n -\gamma_s\Big)\Big|^2ds
	\\ &- n\dint_t^T m_{s-}
	\psi'\Big(m_{s-}(\overline{Y}_{s-}^n -S_{s-})\Big)(\overline{Y}_{s-}^n
	-u_s)^+d\alpha_s - \dint_t^T \psi'\Big(m_{s-}(\overline{Y}_{s-}^n
	-S_{s-})\Big) m_{s-} d\overline{K}_{s}^{n-} \\ &  - \dint_t^T
	\psi'\Big(m_{s}(\overline{Y}_{s}^n -S_{s})\Big) m_{s}\Big(\overline{Z}_{s}^n
	-\gamma_s\Big)dB_s- \dint_t^T\psi'\Big(m_{s-}(\overline{Y}_{s-}^n
	-S_{s-})\Big)\Big(\overline{Y}_{s}^n -S_s\Big)dm_s
	\\ & -\frac12\dint_t^T
	\psi''\Big(m_{s}(\overline{Y}_{s}^n -S_{s})\Big) \Big|m_s\Big(\overline{Z}_{s}^n
	-\gamma_s\Big)\Big|^2ds
	\\ & - \displaystyle{\sum_{t<s\leq T}}
	\bigg[\psi\Big(m_{s}(\overline{Y}_{s}^n
	-S_{s})\Big)-\psi\Big(m_{s-}(\overline{Y}_{s-}^n
	-S_{s-})\Big)-\psi'\Big(m_{s-}(\overline{Y}_{s-}^n -S_{s-})\Big)\Delta_s
	m_.(\overline{Y}_{.}^n -S_.)\bigg].
\end{array}\]
Hence
\begin{equation}\label{equ5}
	\begin{array}{ll}
		\psi\Big(m_t(\overline{Y}_{t}^n -S_t)\Big)&= 1+\dint_t^T
		\psi\Big(m_{s-}(U_{s-} -S_{s-})\Big) m_{s-} dR_s
		%\\&= 1+\dint_t^T
		%\psi(m_{s-}(\overline{Y}_{s-}^n -S_{s-})) m_{s-} dR_s
		%- n\dint_t^T m_{s-}
		%\psi(m_{s-}(\overline{Y}_{s-}^n -S_{s-}))(\overline{Y}_{s-}^n
		%-u_s)^+d\alpha_s\\ & - \dint_t^T
		%\psi(m_{s}(\overline{Y}_{s}^n -S_{s})) m_{s}(\overline{Z}_{s}^n
		%-\gamma_s)dB_s - \dint_t^T \psi(m_{s-}(\overline{Y}_{s}^n
		%-S_{s})) m_{s-} d\overline{K}_{s}^{n-}   \\ &- \dint_t^T\psi(m_{s-}(\overline{Y}_{s-}^n
		%-S_{s-}))(\overline{Y}_{s}^n -S_s)dm_s
		%\\ & - \displaystyle{\sum_{t<s\leq T}}
		%\bigg[\psi(m_{s}(\overline{Y}_{s}^n
		%-S_{s}))-\psi(m_{s-}(\overline{Y}_{s-}^n
		%-S_{s-}))-\psi(m_{s-}(\overline{Y}_{s-}^n -S_{s-}))\Delta_s
		%m_.(\overline{Y}_{.}^n -S_.)\bigg]
		\\ &\quad  -\dint_t^T dV_s^n  - \dint_t^T \psi\Big(m_{s}(\overline{Y}_{s}^n
		-S_{s})\Big) m_{s}\Big(\overline{Z}_{s}^n -\gamma_s\Big)dB_s,
		%\\ &= 1+\dint_t^T
		%\psi(m_{s-}(U_{s-} -S_{s-})) m_{s-} dR_s - n\dint_t^T m_{s-}
		%\psi(m_{s-}(\overline{Y}_{s-}^n -S_{s-}))(\overline{Y}_{s-}^n
		%-u_s)^+d\alpha_s \\ & - \dint_t^T \psi(m_{s}(\overline{Y}_{s-}^n
		%-S_{s})) m_{s}(\overline{Z}_{s}^n -\gamma_s)dB_s-\dint_t^T dV_s^n,
	\end{array}
\end{equation}
where $V^n$ is the process in $\mathcal{K}$ given by
\begin{equation}
\begin{array}{ll}\label{Vn}
	V_t^n = &\dint_0^t
	\Big[\psi\Big(m_{s-}(U_{s-} -S_{s-})\Big)-\psi\Big(m_{s-}(\overline{Y}_{s-}^n -S_{s-})\Big)\Big] m_{s-} dR_s\\ &
	+ n\dint_0^t m_{s-}
	\psi\Big(m_{s-}(\overline{Y}_{s-}^n -S_{s-})\Big)\Big(\overline{Y}_{s-}^n
	-u_s\Big)^+d\alpha_s + \dint_0^t \psi\Big(m_{s-}(\overline{Y}_{s-}^n
	-S_{s-})\Big) m_{s-} d\overline{K}_{s}^{n-}   \\ &+ \dint_0^t\psi\Big(m_{s-}(\overline{Y}_{s-}^n
	-S_{s-})\Big)\Big(\overline{Y}_{s}^n -S_s\Big)dm_s
	\\ & +\sum_{0<s\leq t}\,\psi\Big(m_{s-}(\overline{Y}_{s-}^n -S_{s-})\Big)
	\bigg[\psi\Big(\Delta_s
	m_.(\overline{Y}_{.}^n -S_.)\Big)-1-\Delta_s
	m_.(\overline{Y}_{.}^n -S_.)\bigg].
\end{array}
\end{equation}
Put
\[\begin{array}{ll}
	& \mathcal{M}_t^n = - \Big[\psi\Big(m_t(\overline{Y}_{t}^n -S_t)\Big)+ \dint_0^t
	\psi\Big(m_{u-}(U_{u-} -S_{u-})\Big) m_{u-} dR_u\Big],\\
	& D_t= - \Big[\sup_{s\leqslant t}\psi\Big(m_s(U_{s} -S_s)\Big)+ \dint_0^t
	\psi\Big(m_{u-}(U_{u-} -S_{u-})\Big) m_{u-} dR_u\Big],\\
	& \hat{Z}_{s}^n=\psi\Big(m_{s}(\overline{Y}_{s}^n
	-S_{s})\Big) m_{s}\Big(\overline{Z}_{s}^n -\gamma_s\Big),
\end{array}\]
we obtain from Equation (\ref{equ5}) that
\[ \mathcal{M}_t^n=\mathcal{M}_0^n -V_t^n- \dint_0^t\,\hat{Z}_{s}^ndB_s.\]
Therefore, $(\mathcal{M}_t^{n})_t$  is a \textit{rcll} local
supermartingale satisfying 
\begin{equation}\label{M}
	D_t\leqslant \mathcal{M}_t^n\leqslant \mathcal{M}_t^{n+1}\leqslant0.
\end{equation}
Let $(\tau_i)_{i\geq 0}$ be
the sequence of stopping times defined by
\[\tau_i =\inf\Big\{s\geq 0: -D_s\geq i\Big\}\wedge T.\]
We should note here that the family $(\tau_i)_{i\geq 0}$
satisfies the following property
$$P\bigg[\displaystyle{\bigcup_{i=1}^{\infty}(\tau_i= T)}\bigg]
=1.$$ We say that $(\tau_i)_{i\geq 0}$ is a stationary sequence of
stopping times.
%	\noindent Indeed, let
%	$\displaystyle{\omega\in\bigcap_{j=1}^{\infty}(\tau_j< T)}$ then
%	$\forall j\geq 1, \,\, M_T(\omega):= (A_T+R_T +C_T+\dint_0^T\eta_r
%	dr)(\omega)\geq j$ and hence $M_T(\omega) =+\infty$. Therefore
%	$P\bigg[\displaystyle\bigcap_{j=1}^{\infty}(\tau_j< T)\bigg]\leq
%	P[M_T = +\infty] =0$ and then
%	$P\bigg[\displaystyle{\bigcup_{j=1}^{\infty}(\tau_j= T)}\bigg]
%	=1$.\eop
%\end{remark}
Now, set
\[
\begin{array}{lll}
	& ^i\mathcal{M}_t^n = \mathcal{M}_t^n 1_{\{t<\tau_i\}} + \mathcal{M}_{{\tau_i}-}^n
	1_{\{t\geq\tau_i\}},
	\\ &  ^iV_t^n = V_t^n 1_{\{t<\tau_i\}} + V_{{\tau_i}-}^n
	1_{\{t\geq\tau_i\}}, \\ & ^i\hat{Z}_{s}^n=1_{\{s<\tau_i\}}\hat{Z}_{s}^n .
\end{array}
\]
Then we have
\begin{equation}\label{supermartingale} ^i\!\mathcal{M}_t^n=\mathcal{M}_0^n -^i\!V_t^n- \dint_0^t\,^i\!\hat{Z}_{s}^ndB_s.
\end{equation}
It follows from this last equation that:
\begin{enumerate}
	\item \textbf{$\overline{Y}$ and \underline{Y} are \textit{rcll}}: In fact,
	we have that $(^i\mathcal{M}_t^{n})_t$  is a \textit{rcll} supermartingale satisfying
	$$
	-i\leqslant ^i\mathcal{M}_t^n\leqslant ^i\mathcal{M}_t^{n+1}\leqslant0.
	$$
	It follows then from Dellacherie and Meyer [Theorem 18 Chapiter 5 page 79, \cite{DM}] that
	$\sup_{n}\,\, ^i\mathcal{M}_t^{n}$ is also a
	\textit{rcll} process (supermartingale). Then $\psi(m_t(\overline{Y}_{t}
	-S_t))$ is \textit{rcll} on $[0, \tau_i[$, but
	$(\tau_i)_{i\geq 0}$ is a stationary sequence of stopping times,
	then
	$\overline{Y}$ is \textit{rcll} on $[0, T]$.

	By the same way, we obtain that $\underline{Y}$ is \textit{rcll}.
	\item $\boldmath{
		\overline{Y}_{s-}\leq u_s\,\,\, d\alpha_s- a.e. \,\,\text{on}\,\,
		]0, T]}$: Indeed,  since all terms in equation (\ref{Vn}) are positive it follows that
%	\begin{equation*}\label{equ6}
%	n\int_0^{\tau_i-} m_{s-}
%	\psi\Big(m_{s-}(\overline{Y}_{s-}^n -S_{s-})\Big)\Big(\overline{Y}_{s-}^n
%	-u_s\Big)^+d\alpha_s\leqslant V^n_{\tau_i-}.
%\end{equation*}
%By taking the expectation we have 
\begin{equation}\label{equ6}
n\E\int_0^{\tau_i-} m_{s-}
\psi\Big(m_{s-}(\overline{Y}_{s-}^n -S_{s-})\Big)\Big(\overline{Y}_{s-}^n
-u_s\Big)^+d\alpha_s\leqslant\E V^n_{\tau_i-}.
\end{equation}
%From Equation (\ref{supermartingale}) it follows that 
%$$
%V_{\tau_i-}^n=\mathcal{M}_0^n-\mathcal{M}_{\tau_i-}^n - \dint_0^{\tau_i-}\,\hat{Z}_{s}^ndB_s.
%$$ 
By using a localization procedure and taking the expectation in equation (\ref{supermartingale}) we have 
$$
\E(^i\!V_{T}^n)= \E(\mathcal{M}_0^n-^i\!\mathcal{M}_{T}^n)
$$ 
Hence
$$
\E(V_{\tau_i-}^n)= \E(\mathcal{M}_0^n-\mathcal{M}_{\tau_i-}^n)
$$ 
Since $\mathcal{M}_0^n\leq 0$, it follows from inequality (\ref{M}) and the definition of $\tau_i$ that 
$$
\E V_{\tau_i-}^n= \E(\mathcal{M}_0^n-\mathcal{M}_{\tau_i-}^n) \leq \E(-D_{\tau_i-}) \leq i.
$$ 	
Therefore
\begin{equation*}\label{equ6}
	\E\int_0^{\tau_i-} m_{s-}
	\psi\Big(m_{s-}(\overline{Y}_{s-}^n -S_{s-})\Big)\Big(\overline{Y}_{s-}^n
	-u_s\Big)^+d\alpha_s\leqslant \frac{i}{n}.
\end{equation*}

%	using a standard localization procedure in Equation (\ref{supermartingale}) and the definition of the process $V^n$, we get
	
	Fatou's lemma gives
	\[\E\int_0^{{\tau_{i}}-} m_{s-} \psi\Big(m_{s-}(\overline{Y}_{s}^{\,\,-}
	-S_{s-})\Big)\Big(\overline{Y}_{s}^{\,\,-} -u_s\Big)^+d\alpha_s =0.\]
	Hence
	\[\overline{Y}_{s}^{\,\,-}\leq u_s\,\,\, d\alpha_s- a.e.
	\,\,\text{on}\,\, ]0, T[.\]
	But for every $s\in ]0, T]$ and $n\in\N$,\,\,$\overline{Y}_{s-}\leq
	\overline{Y}_{s-}^n$ then $\overline{Y}_{s-}\leq
	\overline{Y}_{s}^{\,\,-}$. Consequently
	$$
	\overline{Y}_{s-}\leq u_s\,\,\, d\alpha_s- a.e. \,\,\text{on}\,\,
	]0, T[.
	$$
	Assume now that $\overline{Y}_{T-} > u_T$ and $\Delta_T\alpha >0$.
	It follows from [\cite{EH4}, Lemma 3.1.] that
	$$
	\overline{Y}_{T-}^n = S_{T^-}\vee \Big[\xi+\Delta_TV^++\Delta_TV^-+\beta_T\Delta_T A
	-n\Big(\overline{Y}_{T-}^n -u_T\Big)^+\Delta_T\alpha\Big]\wedge U_{T-}.
	$$
	Since $-n(\overline{Y}_{T-}^n -u_T)^+$ converges to $-\infty$ if $n$
	goes to $+\infty$, we get $\overline{Y}_{T}^{\,\,-} = S_{T-}$.
	Now since $S_{T-}\leq \overline{Y}_{T-}\leq
	\overline{Y}_{T}^{\,\,-} $ we have $\overline{Y}_{T-}= S_{T-}$.
	Hence $S_{T-}> u_T$ which is absurd since $\Delta_T\alpha >0$ and
	$S_{t-}\leq u_t$ $d\alpha_t- a.e. \,\,\text{on}\,\, ]0, T].$
	Consequently
	$$
	\overline{Y}_{s-}\leq u_s\,\,\, d\alpha_s- a.e. \,\,\text{on}\,\,
	]0, T].
	$$
	\item By the same method as in the previous step, we get also that
	$$
	l_t\leq \underline{Y}_{t-}, \,\,\, d\delta_t-a.e., \text{on}\,\, ]0,T].
	$$ Hence $\overline{Y}$ and $\underline{Y}$ are in $Dom$. The proof of Proposition \ref{pro1} is finished.
\end{enumerate}
\eop
\subsection{Proof of Theorem \ref{thee1}}
With the help of processes $\overline{Y}$ and $\underline{Y}$, another RBSDE is considered which is equivalent to our original RBSDE. More precisely, we have the following proposition.
\begin{proposition}
	$(Y, Z, K^+, K^-)$ is a solution of RBSDE (\ref{eq0}) if and only
	if it is a solution of the following  RBSDE%(\ref{eq0})
	\begin{equation}
		\label{eq01} \left\{
		\begin{array}{ll}
			(i) &  Y_{t}=\xi
			+\integ{t}{T}\left[f(s,Y_{s},Z_{s})ds+g(s,
			Y_{s-}, Y_s)dA_s+dK_{s}^+ -dK_{s}^-
			-Z_{s}dB_{s}\right], t\leq T,
			\\ (ii)&
			\forall t\in[0,T[,\,\, \underline{Y}_t \leq Y_{t}\leq
			\overline{Y}_{t},
			\\ (iii)&   \integ{0}{T}(
			Y_{t-}-\underline{Y}_{t-})
			dK_{t}^+= \integ{0}{T}( \overline{Y}_{t-}-Y_{t-}) dK_{t}^-=0,\,\, \text{a.s.}, \\
			(v)& Y\in {\cal D}, \quad K^+, K^-\in {\cal K}, \quad Z\in {\cal
				L}^{2,d},  \\ (vi)& dK^+\perp  dK^-.
		\end{array}
		\right. \end{equation}
\end{proposition}

\bop. Let $(Y, Z, K^+, K^-)$ be a solution of RBSDE (\ref{eq0}). Put
$$
U_t^* = Y_t\vee \overline{Y}_t\quad \text{and} \quad L_t^* =
Y_t\wedge \underline{Y}_t.
$$
It is obvious that $L_t^*\leq Y_t\leq U_t^*$. By Remark \ref{rem20}, it follows that $U^*, L^*$  are in $Dom$. Then $(Y, Z, K^+, K^-)$ is a solution
of the following RBSDE
\begin{equation*}
	\label{eq0100} \left\{
	\begin{array}{ll}
		(i) & Y_{t}=\xi
		+\integ{t}{T}\left[f(s,Y_{s},Z_{s})ds+g(s,
		Y_{s-}, Y_s)dA_s+dK_{s}^+ -dK_{s}^-
		-Z_{s}dB_{s}\right], t\leq T,
		\\ (ii)&
		\forall t\in[0,T[,\,\, L^*_t \leq Y_{t}\leq U_{t}^*,
		\\ (iii)&   \integ{0}{T}(
		Y_{t-}-L_{t-}^*)
		dK_{t}^+= \integ{0}{T}( U_{t-}^*-Y_{t-}) dK_{t}^-=0,\,\, \text{a.s.}, \\
		(v)& Y\in {\cal D}, \quad K^+, K^-\in {\cal K}, \quad Z\in {\cal
			L}^{2,d},  \\ (vi)& dK^+\perp  dK^-.
	\end{array}
	\right. \end{equation*} Since
\begin{itemize}
	\item[(a)] $Y_s\leqslant U_s$ and $L_s^*\leqslant\underline{Y}_t\leqslant  S_t\leqslant \overline{Y}_s^n$,
	\item[(b)]
	$$
	\begin{array}{lll}
	& f(s,Y_{s},Z_{s})ds+g(s, Y_{s-}, Y_s)dA_s \\ &\leqslant (\eta_s +C_s|Z_s|^2 )ds+\beta_sdA_s
 \\ & \leqslant (\eta_s
	+4C_s|\gamma_s|^2 + \frac{m_s}{2}|Z_s-\gamma_s|^2)ds +\beta_sdA_s +dV^+_s+dV^-_s-\underbrace{n(u_s-Y_{s-})^-d\alpha_s}_{=0},\end{array}
$$
\end{itemize}

then it follows from comparison theorem (Theorem \ref{th2}) applied to $Y$ and $\overline{Y}^n$, that for all $n\in\N$, $$Y_t\leq
\overline{Y}_t^n,$$
and then $Y_t\leq \overline{Y}_t$. Hence $U_t^*
=\overline{Y}_t$. By  a symmetric argument we get also that $L_t^*
=\underline{Y}_t$. Therefore $(Y, Z, K^+, K^-)$ is a solution to
RBSDE (\ref{eq01}).

Conversely, suppose now that $(Y, Z, K^+, K^-)$ is a solution of
(\ref{eq01}).  In
order to prove that $(Y, Z, K^+, K^-)$ is a solution of RBSDE
(\ref{eq0}) we just need to prove $(iii)$  of RBSDE (\ref{eq0}). Let $L^*\in Dom$ and consider
$(Y^*, Z^*, K^{+*}, K^{-*})$ the minimal solution of the following RBSDE
\begin{equation*}
	\label{eq01000} \left\{
	\begin{array}{ll}
		(i) &  Y_{t}^*=\xi
		+\integ{t}{T}[(f(s,Y_{s},Z_{s})-\frac{m_s}{2}|Z_{s}^*-Z_s|^2)ds+g(s,
		Y_{s-}, Y_s)dA_s +dK_{s}^{*+} -dK_{s}^{*-}
		-Z_{s}^*dB_{s}],
		\\ (ii)&
		\forall t\in[0,T[,\,\, {Y}_t \leq Y_{t}^*\leq L^*_t\vee Y_t,
		\\ (iii)&   \integ{0}{T} \Big(
		Y_{t-}^*-Y_{t-} \Big)
		dK_{t}^{*+}= \integ{0}{T} \Big((L_{t-}^*\vee Y_{t-})-Y_{t-}^*\Big) dK_{t}^{*-}=0,\,\, \text{a.s.}, \\
		(v)& Y^*\in {\cal D}, \quad K^{*+}, K^{*-}\in {\cal K}, \quad Z^*\in
		{\cal L}^{2,d},  \\ (vi)& dK^{*+}\perp  dK^{*-}.
	\end{array}
	\right. \end{equation*}
We note here that this minimal solution exists according to \cite{EH4}. On the other hand we have

\begin{itemize}
	\item[(a)] $Y_s\leq
	\overline{Y}_s\leq\overline{Y}_s^n $ and $Y^*_s\leq L^*_s\vee Y_s\leq L^*_s\vee \overline{Y}_s\leq U_s$,
	\item[(b)] $[f(s,Y_{s},Z_{s})-\frac{m_s}{2}|Z_{s}^*-Z_s|^2]ds+g(s, Y_{s-}, Y_s)dA_s \leqslant (\eta_s +C_s|Z_s|^2 -\frac{m_s}{2}|Z_{s}^*-Z_s|^2)ds+\beta_sdA_s$
	\item[] \qquad$\leqslant (\eta_s
	+4C_s|\gamma_s|^2 + \frac{m_s}{2}|Z^*_s-\gamma_s|^2)ds +\beta_sdA_s +dV^+_s+dV^-_s-\underbrace{n(u_s-Y^*_{s-})^-d\alpha_s}_{=0}$.
\end{itemize}
Applying the comparison theorem (Theorem \ref{th2})  to $Y^*$ and
$\overline{Y}^n$, we get $\forall n\in\N$,\,\,$Y^*\leq
\overline{Y}^n$. Letting $n$ to infinity we obtain
$
Y^*\leq \overline{Y}.
$
Applying again comparison theorem (Theorem \ref{th2}) to $-Y^*$ and $-Y$ it follows that
$Y^*\leq Y$.
Then
\[Y^*= Y, \,\,\, Z= Z^*\,\text{and}\,\,\,dK^{*-} = dK^-.\]
Henceforth
\[(Y_{s-}-L_{s-}^*)^-\,dK^-_s =(L_{s-}^* \vee Y_{s-}-Y_{s-})dK^{*-}_s = 0.\]
By symmetric argument we get also,
\[(Y_{s-}-L_{s-}^*)^+dK^+_s = 0.\]
Consequently $(Y, Z, K^+, K^-)$ is a solution of (\ref{eq0}).\eop

As consequence we get the following result
\begin{corollary}\label{pro2}
	%$(Y, Z, K^+, K^-)$ is a solution of RBSDE (\ref{eq0}) if and only if $(Y, Z, K^{+}, K^{-})$ is a solution of RBSDE (\ref{eq01}). Hence
	$(Y, Z, K^+, K^-)$ is a maximal (resp. minimal) solution of RBSDE (\ref{eq0}) if and only if $(Y, Z, K^{+}, K^{-})$ is a maximal (resp. minimal) solution of RBSDE (\ref{eq01}).
\end{corollary}
\bop\,\,\textbf{of Theorem \ref{thee1}}. According to [Theorem 2.1., \cite{EH4}], there exists $(Y, Z, K^+, K^-)$  a maximal (resp. minimal) solution of RBSDE (\ref{eq01}) and then by Corollary \ref{pro2}, $(Y, Z, K^+, K^-)$ is a maximal (resp. minimal) solution of RBSDE
(\ref{eq0}).\eop

\section{Further study: standard form of RBSDE}\label{sec2}
In this section, we want to find an equivalent and standard form to our initial RBSDE (\ref{eq0}) by giving another characterization of the $Dom$ without introducing the test barriers $L^*$ and $U^*$.\\

First, from the appendix, we have the following characterization of $Dom$.
\begin{proposition}
	$Dom=\bigg\{Y\in \mathcal{D}\,\,:\;\; \left[\forall t\in ]0,T],\; L_{t-}\vee l^{*,\delta}_t \leqslant Y_{t-}\leqslant U_{t-}\wedge (-(-u)^{*,\alpha}_t)\right]\; P\text{-as}\bigg\}$,
	where $l^{*,\delta}(t)$ and $(-u)^{*,\alpha}(t)$ are defined respectively by (see the appendix for more details)
	$$
	l^{*,\delta}(t)=	\inf_{n}\bigg\{-nt+ \inf\Big\{ a\in \R: \dint_0^t\Big[l_s+ns-
	a\Big]^+\,d\delta_s=0\Big\}\bigg\},
	$$
	$$
	(-u)^{*,\alpha}(t)=	\inf_{n}\bigg\{-nt+ \inf\Big\{ a\in \R: \dint_0^t\Big[-u_s+ns-
	a\Big]^+\,d\alpha_s=0\Big\}\bigg\}.
	$$
\end{proposition}

The following theorem  proves that our original RBSDE can be written in a standard form.
\begin{theorem}
	$(Y, Z, K^+, K^-)$ is a solution of RBSDE (\ref{eq0}) if and only if  $(Y, Z, K^+, K^-)$ satisfies
	\begin{equation}
		\label{eq*} \left\{
		\begin{array}{ll}
			(i) &  Y_{t}=\xi
			+\integ{t}{T}\left[f(s,Y_{s},Z_{s})ds+g(s,
			Y_{s-}, Y_s)dA_s+dK_{s}^+ -dK_{s}^-
			-Z_{s}dB_{s}\right],
			\\ (ii)&
			\forall t\in]0,T],\,\, L_{t-}\vee l^{*,\delta}(t)\leqslant Y_{t-}\leqslant -(-u)^{*,\alpha}(t)\wedge U_{t-},\\  (iii)&
			\integ{0}{T}( Y_{t-}-[L_{t-}\vee l^{*,\delta}(t)])
			dK_{t}^+= \integ{0}{T}( [U_{t-}\wedge -(-u)^{*,\alpha}(t)]-Y_{t-}) dK_{t}^-=0,\,\, \text{a.s.}, \\
			(iv)& Y\in {\cal D}, \quad K^+, K^-\in {\cal K}, \quad Z\in {\cal
				L}^{2,d},  \\ (v)& dK^+\perp  dK^-.
		\end{array}
		\right. \end{equation}
\end{theorem}
\bop.

Let $(Y, Z, K^+, K^-)$ be a solution of RBSDE (\ref{eq0}), then for all processes $ L^*$ and $U^*$ in $Dom$ such that for all $t\in [0,T[$, $L^*_t\leqslant Y_t\leqslant U^*_t$ we have $(Y, Z, K^+, K^-)$ is also a solution of the following RBSDE
\begin{equation*}
	\label{eq**} \left\{
	\begin{array}{ll}
		(i) &  Y_{t}=\xi
		+\integ{t}{T}\left[f(s,Y_{s},Z_{s})ds+g(s,
		Y_{s-}, Y_s)dA_s+dK_{s}^+ -dK_{s}^-
		-Z_{s}dB_{s}\right],
		\\ (ii)&
		\forall t\in[0,T[,\,\, L^*_t \leq Y_{t}\leq U^*_{t},\\  (iii)&
		\integ{0}{T}( Y_{t-}-L^*_{t-})
		dK_{t}^+= \integ{0}{T}( U^*_{t-}-Y_{t-}) dK_{t}^-=0,\,\, \text{a.s.}, \\
		(iv)& Y\in {\cal D}, \quad K^+, K^-\in {\cal K}, \quad Z\in {\cal
			L}^{2,d},  \\ (v)& dK^+\perp  dK^-.
	\end{array}
	\right. \end{equation*}

But for each integers $n$ and $m$,
$$
L^n_t(\omega)=L_t(\omega)\vee l^{n,\delta}(t+,\omega)\wedge Y_t(\omega)\quad\text{and} \quad U^m_t(\omega)=Y_t(\omega)\vee [-(-u)^{m,\alpha}(t+,\omega)]\wedge U_t(\omega),$$
are in $Dom$ and $L^n_t\leqslant Y_t\leqslant U^m_t$, where 
$$
l^{n,\delta}(t,\omega) = -nt+\essup^{\delta}_{s\leqslant t}\, \Big[l(s,\omega)+ns\Big]\,\, \text{ and }\,\, (-u)^{m,\alpha}(t,\omega) = -mt+\essup^{\alpha}_{s\leqslant t}\, \Big[-u(s,\omega)+ms\Big].
$$
%(See appendix for properties of the processes $l^{n,\delta}$ and $(-u)^{m,\alpha}$).

We deduce that, for each integers $n$ and $m$
\begin{align*}
	&\integ{0}{T}( Y_{t-}-L^n_{t-})dK_{t}^+= \integ{0}{T}( U^m_{t-}-Y_{t-}) dK_{t}^-=0.
	%& Y_{t-}=L^n_{t-}\vee \big(Y_t+g(t,Y_{t-},Y_t)\Delta_tA\big)\wedge U^m_{t-},\\
	%&\Delta_tK^+ =\bigg(L^n_{t-}- \big(Y_t+g(t,Y_{t-},Y_t)\Delta_tA\big)\bigg)^+,\\
	%&\Delta_tK^- =\bigg(U^m_{t-}- \big(Y_t+g(t,Y_{t-},Y_t)\Delta_tA\big)\bigg)^-.
\end{align*}
It follows from Appendix that
$$
l^{*,\delta}_-(t)\leqslant l^{*,\delta}(t)\leqslant Y_{t-}\leqslant -(-u)^{*,\alpha}(t)\leqslant -(-u)_-^{*,\alpha}(t).
$$
Then passing in limit and using monotone convergence theorem we get
\begin{align*}
	&\integ{0}{T}\Big(Y_{t-}-[L_{t-}\vee l^{*,\delta}_-(t)]\Big)dK_{t}^+= \integ{0}{T}\Big((U_{t-}\wedge -(-u)_-^{*,\alpha}(t)) -Y_{t-}\Big) dK_{t}^-=0.
	%& Y_{t-}=[L_{t-}\vee l^{*,\delta}_-(t)]\vee \big(Y_t+g(t,Y_{t-},Y_t)\Delta_tA\big)\wedge [U_{t-}\wedge -(-u)_-^{*,\alpha}(t)],\\
	%&\Delta_tK^+ =\bigg([L_{t-}\vee l^{*,\delta}_-(t)]- \big(Y_t+g(t,Y_{t-},Y_t)\Delta_tA\big)\bigg)^+,\\
	%&\Delta_tK^- =\bigg([U_{t-}\wedge -(-u)_-^{*,\alpha}(t)]- \big(Y_t+g(t,Y_{t-},Y_t)\Delta_tA\big)\bigg)^-.
\end{align*}
This gives the necessary implication. \\
 Let us now show the reverse. Suppose that $(Y, Z, K^+, K^-)$ is a solution of RBSDE (\ref{eq*}). In order to prove that $(Y, Z, K^+, K^-)$ is a solution of RBSDE (\ref{eq0})  it remains to prove $(iii)$ and $Y\in Dom$. Let $L^*\in Dom$, it follows that
$$
\bigg(l_t(\omega)\leqslant L^*_{t-}(\omega)\bigg) \;\; d\delta_t(\omega)P(d\omega)\text{-a.e} \,\, \text{and}\,\, L_{t-} \leqslant L^*_{t-}\,\,dt\,P(d\omega)\text{-a.e}.
$$
Since $L$ and $L^*$  are \it{rcll}, by using Corollary  \ref{char}, we have
%the characterization: $\forall Y\in Dom$
%\[\bigg(l_t(\omega)\leqslant Y_{t-}(\omega)\bigg) \;\; d\delta_t(\omega)P(d\omega)\text{-a.e} \quad\text{ if and only if }\quad \bigg(\forall t\in ]0,T],\;\; l_t^{*,\delta}(\omega)\leqslant Y_{t-}(\omega)\bigg)\;\; P\text{-as} \]
%it follows that 
$$
\bigg(\forall t\in ]0,T],\;\; l_t^{*,\delta}\leqslant L^*_{t-} \,\, \text{and}\,\, L_{t-} \leqslant L^*_{t-}\bigg)\,\,\;\; P\text{-as}
$$
Henceforth 
$$
\bigg(\forall t\in ]0,T],\;\; l_t^{*,\delta}\vee L_{t-}\leqslant L^*_{t-}\bigg)\,\,\;\; P\text{-as}
$$
Using the same method as above we have also 
$$
\bigg(\forall t\in ]0,T],\;\;  L^*_{t-}\leqslant -(-u)^{*,\alpha}(t)\wedge U_{t-}\bigg)\,\,\;\; P\text{-as}
$$
Henceforth
$$
\begin{array}{ll}
	& \E\dint_0^T\;(
	Y_{t-}-L_{t-}^*)^+\,dK_{t}^++(
	Y_{t-}-L_{t-}^*)^-\,dK_{t}^- \\ & \leq \E\integ{0}{T}( Y_{t-}-[L_{t-}\vee l^{*,\delta}(t)])
	dK_{t}^+ +\E\integ{0}{T}( [U_{t-}\wedge -(-u)^{*,\alpha}(t)]-Y_{t-}) dK_{t}^- =0.
\end{array}
$$
Then $(iii)$ is satisfied.

 Let us now show that $Y\in Dom$.
It follows from $(ii)$ and Corollary \ref{char} that
$$
\E\dint_0^T\left(L_{t-}-Y_{t-}\right)^+dt+\left(l_{t}-Y_{t-}\right)^+d\delta_t=0.
$$
Now by taking $-Y$ and $-u$, in Corollary \ref{char}, instead of $Y$ and $l$ we have also 
$$
\E\dint_0^T\left(U_{t-}-Y_{t-}\right)^-dt+\left(u_{t}-Y_{t-}\right)^-d\alpha_t=0.
$$
This gives the result.
\eop
\begin{remark}
	From the above proof, we have
	\[\integ{0}{T}\Big([L_{t-}\vee l^{*,\delta}(t)]-[L_{t-}\vee l^{*,\delta}_-(t)]\Big)dK_{t}^+= \integ{0}{T}\Big([U_{t-}\wedge -(-u)_-^{*,\alpha}(t)]- [U_{t-}\wedge -(-u)^{*,\alpha}(t)]\Big) dK_{t}^-=0\]
	which is equivalent to
	\[
	\text{if}\,\, l_t>[L_{t-}\vee l^{*,\delta}_-(t)] \text{ and } \Delta_t\delta>0\,\,  \text{then}\,\,\Delta_tK^+= 0,\]and
	\[
	\text{if}\,\,  u_t<[U_{t-}\wedge -(-u)_-^{*,\alpha}(t)]\text{ and } \Delta_t\alpha>0\,\,  \text{then}\,\,\Delta_tK^-= 0.\]
	%\[\forall t\in]0,T],\;\;([L_{t-}\vee l^{*,\delta}(t)]-[L_{t-}\vee l^{*,\delta}_-(t)])\Delta_tK^+= ([U_{t-}\wedge -(-u)_-^{*,\alpha}(t)] -[U_{t-}\wedge -(-u)^{*,\alpha}(t)]) \Delta_tK^-=0\]
\end{remark}

\section{Particular case: Generalized Snell envelope}\label{sec}

Let $l$\ be a predictable process, $L\in {\cal D}$  and $\delta\in {\cal K}$ satisfying the following
hypothesis:

\noindent \textbf{(A)} There exists a local martingale $M_t = M_0 + \dint_0^t
\kappa_s dB_s$ such that $P-$a.s.,  \[L_t\leq M_t \;\text{ on }\;
]0, T]\quad\text{and} \quad l_t\leq M_{t}
\;d\delta_t\text{-a.e.\;on } ]0, T].\]

\noindent Let $U_t = U_0 -V_t +\dint_0^t \chi_s dB_s$, where $V\in {\cal
	K}$ and $\chi \in {\cal L}^{2, d}$, be a \textit{rcll} local
supermartingale such that
$P-$a.s.,  \[L_t\leq U_t \;\text{ on }\;
]0, T]\quad\text{and} \quad l_t\leq U_{t-}
\;d\delta_t\text{-a.e.\;on } ]0, T].\]
According to our main result, let $(Y, Z, K^+, K^-)$  be the minimal solution of the following
RBSDE
\begin{equation*}
	\label{eq00100} \left\{
	\begin{array}{ll}
		(i) &  Y_{t}=L_T+\integ{t}{T}\Big[dK_{s}^+ -dK_{s}^--Z_{s}dB_{s}\Big],
		\\ (ii)&
		\forall t\in]0,T],\,\, L_{t-}\vee l^{*,\delta}(t)\leqslant Y_{t-}\leqslant  U_{t-},\\  (iii)&
		\integ{0}{T}( Y_{t-}-[L_{t-}\vee l^{*,\delta}(t)])
		dK_{t}^+= \integ{0}{T}( U_{t-}-Y_{t-}) dK_{t}^-=0,\,\, \text{a.s.}, \\
		(iv)& Y\in {\cal D}, \quad K^+, K^-\in {\cal K}, \quad Z\in {\cal
			L}^{2,d},  \\ (v)& dK^+\perp  dK^-.
	\end{array}
	\right. \end{equation*}
Applying comparison theorem (Theorem \ref{th2})  to the processes $U_t$ and $Y_t$, we get $dK^- =0$. Then $Y$ is a local supermartingale
minimal solution of the following RBSDE
\begin{equation}
	\label{eq001000} \left\{
	\begin{array}{ll}
		(i) &  Y_{t}=L_T+\integ{t}{T}\Big[dK_{s}^+ -Z_{s}dB_{s}\Big],
		\\ (ii)&
		\forall t\in]0,T],\,\, L_{t-}\vee l^{*,\delta}(t)\leqslant Y_{t-},\\  (iii)&
		\integ{0}{T}( Y_{t-}-[L_{t-}\vee l^{*,\delta}(t)])
		dK_{t}^+=0,\,\, \text{a.s.}, \\
		(iv)& Y\in {\cal D}, \quad K^+\in {\cal K}, \quad Z\in {\cal
			L}^{2,d}.
	\end{array}
	\right. \end{equation}

Henceforth, we obtain the following result.
\begin{theorem} Suppose that \textbf{(A)} hold. Then
	 the minimal solution $Y$ of (\ref{eq001000}) is the smallest
	\textit{rcll} local super-martingale satisfying
	$P-$a.s.,  \[L_t\leq Y_t \;\text{ on }\;
	]0, T]\quad\text{and} \quad l_t\leq Y_{t-}
	\;d\delta_t\text{-a.e.\;on } ]0, T].\]
\end{theorem}
We say that $Y$ is the generalized Snell envelope associated to $L$ and $ld\delta$. We denote it by $\mathcal{S}_.(L, ld\delta)$.
\begin{remark} We know that if $L$ is of class $\it D$ then $L$ satisfies assumption
	\textbf{(A)} (see Dellacherie-Meyer \cite{DM}, Theorem 24 page 419). In this case our generalized Snell envelope
	$\mathcal{S}_.(L) = \mathcal{S}_.(L, 0d0)$ coincides with the
	usual Snell envelope $\essup_{\tau \in {\cal T}_{t}}
	\E[L_{\tau}|{\cal F}_t]$, where ${\cal T}_{t}$ is the set of all
	stopping times valued between $t$ and $T$, as presented in
	[Dellacherie-Meyer \cite{DM}, page 416] and studied by several authors.
\end{remark}
Let us give the following two examples in the case where $\delta_t =\lambda$ the Lebesgue measure and $\delta_t = 1_{\{T'\leq t\}}$, where $T'$ is a stopping time with values in $[0,T]$.
\begin{example} \noindent Let $l$\ be a predictable process and $\xi$ a $\mathcal{F}_T$-measurable random variable such that there exist $L\in {\cal D}$ and $M$ a local martingale such that
	$L_{t}\leq l_t\leq M_t\; d\lambda\text{-a.e}$ and $\xi\leq M_T$ (where $\lambda$ denotes the Lebesgue measure). Let $(Y, Z, K^+, K^-)$ be
	the minimal solution of the following RBSDE
	$$
	\label{eq010000} \left\{
	\begin{array}{ll}
		(i) & %\text{ is a f-supersoltion on } [0, T]\text{ with }Y_T = \xi, i.e.
		Y_{t}=\xi
		+\dint_{t}^{T}\Big[dK_{s}^+ -Z_{s}dB_{s}\Big],\; t\leq T,
		\\ (ii)& l^{*,\lambda}_t \leq Y_{t-}\;\text{ on }\,\, ]0, T],
		\\ (iii)&   \E\dint_{0}^{T}( Y_{t-}-l^{*,\lambda}_t)
		dK_{t}^+=0, \\
		(iv)& Y\in {\cal D}, \quad K^+\in {\cal K}, \quad Z\in {\cal
			L}^{2,d},
	\end{array}
	\right. $$ Then $Y$ $\big(=\mathcal{S}_.(L_t1_{\{t<T\}}+\xi1_{\{t=T\}}, ld\lambda)\big)$ is the smallest local super-martingale such that
	$$
	l_t\leq Y_t,\,\,d\lambda- a.e\,\,\text{ and }\,\, \xi\leq Y_T.
	$$
   % Note that for every $t\in]0,T]$, $L_{t-}\leqslant l^{*,\lambda}_t$.
\end{example}
\begin{example} \noindent Let $T'$ be a stopping time with values in $[0,T]$, $\xi'$ a $\mathcal{F}_{T'}$-measurable random variable,  $\xi$ a $\mathcal{F}_T$-measurable random variable and $L\in {\cal D}$ such that there exists  a local martingale $M$ such that
	$L_{t}\leq M_t\; \text{ on } [0,T[$, $\xi'\leq M_{T'}$ and $\xi\leq M_T$. Let $(Y, Z, K^+, K^-)$ be
	the minimal solution of the following RBSDE
	$$
	\label{eq010000x} \left\{
	\begin{array}{ll}
		(i) & Y_{t}=\xi+\dint_{t}^{T}\Big[dK_{s}^+ -Z_{s}dB_{s}\Big],\; t\leq T,
		\\ (ii)& L_{t-}\vee l^{*,\delta}_t \leq Y_{t-}\;\text{ on }\,\, ]0, T],
		\\ (iii)&   \E\dint_{0}^{T}( Y_{t-}-L_{t-}\vee l^{*,\delta}_t)
		dK_{t}^+=0, \\
		(iv)& Y\in {\cal D}, \quad K^+\in {\cal K}, \quad Z\in {\cal
			L}^{2,d},
	\end{array}
	\right. $$ where $\delta_t=1_{\{T'\leqslant t\}}$ and $l_t=\xi'1_{\{t=T'\}}$ hence $ l^{*,\delta}_t=-\infty\mathbf{1}_{\{t\neq T'\}}+\xi'\mathbf{1}_{\{t=T'\}}$ . Then $Y$ ($=\mathcal{S}_.(L_t1_{\{t<T\}}+\xi1_{\{t=T\}}, ld\delta)$) is the smallest local super-martingale such that
	\[L_t\leq Y_t,\;\text{ on }[0,T[\;,\quad  \xi'\leq Y_{T'-}\quad\text{ and }\quad \xi\leq Y_T.\]
\end{example}

%\begin{remark}
%The difference between the above definition and those of continuous or \it{rcll} obstacles or even irregular $L^2$-obstacles, is, on one hand, the dominating conditions $(ii)-(iii)$, and on another hand the Skorohod condition $(iv)$.
%\end{remark}

\section{\Large{\textbf{Appendix}}}

In this appendix we give, in particular, the following characterization: for  $Y\in\mathcal{D}$ we have the following
\[\big(g(t,\omega)\leqslant Y_{t-}(\omega)\big) \;\; d\rho_t(\omega)P(d\omega)\text{-a.e} \quad\text{ if and only if }\quad \big(\forall t\in ]0,T],\;\; g^{*,\rho}(t,\omega)\leqslant Y_{t-}(\omega)\big)\;\; P\text{-as,} \]
where $\rho$ is a process in $\mathcal{K}$.
This characterization allow us to write our RBSDE in standard form.  

Let $g: ]0, T]\times \Omega \rightarrow \R$ be a progressively
measurable function and $\rho$ be a process in $\mathcal{K}$.
%such that for every $\omega$,
%$d\delta_t(\omega)-{\essup^\delta_{t\in ]0, T]}} g(t,\omega)
%:= \Psi(\omega) < +\infty$.

Note for $(t,\omega)\in[0, T]\times\Omega$ and $n\in\N$
\begin{equation}\label{apen}
	\begin{array}{lll}	
g^{n,\rho}(t,\omega)& =-nt+\displaystyle{\essup^\rho_{s\leqslant t}}\, \Big[g(s,\omega)+ns\Big] \\ & =-nt+ \inf\bigg\{ \alpha\in \R: \dint_0^t\Big[g(s,\omega)+ns-
\alpha\Big]^+\,d\rho_s(\omega)=0\bigg\}.
\end{array}
\end{equation}
We have the following.
\begin{proposition} \label{prop1}For every $n\in \N^*$, $g^{n,\rho}$ is a predictable process satisfying:
	\begin{enumerate}
		\item For each $\omega\in\Omega$, \,$t\mapsto nt+g^{n,\rho}(t,\omega)$ is a non-decreasing function and then $t\mapsto g^{n,\rho}(t,\omega)$ is làglàd with $g^{n,\rho}(t-,\omega)\leqslant g^{n,\rho}(t,\omega) \leqslant g^{n,\rho}(t+,\omega)$ and such that
		\begin{enumerate}
			\item $g^{n,\rho}(t-,\omega)= -nt+\inf\Big\{ \alpha\in \R: \dint_0^{t-}\Big[g(s,\omega)+ns-
			\alpha\Big]^+\,d\rho_s(\omega)=0\Big\}.$
			\item If $\Delta_t\rho(\omega)>0$, then $g^{n,\rho}(t,\omega)=g^{n,\rho}(t-,\omega)\vee g(t,\omega)$.
			\item If $\Delta_t\rho(\omega)=0$, then $g^{n,\rho}(t,\omega)=g^{n,\rho}(t-,\omega).$
		\end{enumerate}
		\item For each $(t,\omega)\in[0, T]\times\Omega$, \,$-\infty\leq g^{n+1,\rho}(t,\omega)\leq g^{n,\rho}(t,\omega)$.
		\item If $g\leqslant h\,\;\,d\rho_s(\omega)P(d\omega)$--a.e. then $P$-a.e. $\omega\in\Omega$  for every $t\in[0, T]$ and $n\in\N$, we have \[g^{n,\rho}(t,\omega)\leqslant h^{n,\rho}(t,\omega)\leqslant \sup_{s\leqslant t}\, \Big[h(s,\omega)-n(t-s)\Big].\]
	\end{enumerate}
\end{proposition}
\bop. Properties 1, 2 and 3 are obvious. The predictability follows from the fact that
\[ \Big\{(t,\omega): nt+g^{n,\rho}(t,\omega)\leq a\Big\} = \Big\{(t,\omega): \dint_0^t \Big[g(s,\omega) + ns -a\Big]^+ d\rho_s(\omega) =0\Big\}.\]\eop
Let us now define
$$
g^{*,\rho}(t,\omega) = \inf_n
g^{n,\rho}(t,\omega),\quad \,g_-^{*,\rho}(t,\omega) = \inf_n
g^{n,\rho}(t-,\omega)\quad \text{and}\quad g_+^{*,\rho}(t,\omega) = \inf_n
g^{n,\rho}(t+,\omega).$$
\begin{remark}
	We have the following:
	\begin{enumerate}
		\item  $g_-^{*,\rho}(t,\omega)\leqslant g^{*,\rho}(t,\omega) \leqslant g_+^{*,\rho}(t,\omega)$.
		\item  If $\Delta_t\rho(\omega)>0$, then  $g^{*,\rho}(t,\omega)=g_-^{*,\rho}(t,\omega)\vee g(t,\omega)$.
		\item  If $\Delta_t\rho(\omega)=0$ then $g^{*,\rho}(t,\omega)=g_-^{*,\rho}(t,\omega).$
	\end{enumerate} In particular, for every $\omega\in \Omega$ \[\Big\{t\in]0,T]\,\,: \,\, g^{*,\rho}(t,\omega)>g_-^{*,\rho}(t,\omega)\Big\}\subset \Big\{t\in]0,T]\,\,: \,\, \Delta_t\rho(\omega)>0\Big\}\] and then
	$\{t\in]0,T]\,\,: \,\, g^{*,\rho}(t,\omega)>g_-^{*,\rho}(t,\omega)\}$ is a countable set.
\end{remark}

\begin{proposition} \label{prop3} For every $\omega\in\Omega$,
	$$
	g(t, \omega) \leq g^{*,\rho}(t,\omega) \,\,\,d\rho_t(\omega)-a.e.\quad \text{on}\quad ]0, T].
	$$
\end{proposition}
\bop. From the definition of $g^{n,\rho}$, we have that for every $(\omega,t ,n)$ in $\Omega\times]0,T]\times\N$, there exists a negligible Borel
set $N_{(\omega,t ,n)}$ with respect to the measure $d\rho_s(\omega)$ such that for every $s\in N_{(\omega,t ,n)}^c$
$$
1_{\{s\leq t\}} g(s,\omega)\leq \Big[g^{n,\rho}(t,\omega)+n(t-s)\Big]1_{\{s\leq t\}}.
$$
For $\omega\in\Omega$, let
$$
I_{\omega} = \bigcup_{n\in\N}
\Big\{ t\in ]0, T[ \,\,:\,\, g^{n,\rho}(t-,\omega) < g^{n,\rho}(t+,\omega)\Big\}\bigcup \bigg[\Q\bigcap ]0, T[\bigg]\bigcup \Big\{T\Big\},
$$
which is countable and dense in $[0, T]$.\\ Define the following negligible Borel set  with respect to the measure $d\rho_s(\omega)$
$$
N_{\omega} = \bigcup_{t\in I_{\omega}}\bigcup_{n\in\N} N_{(\omega,t ,n)}.
$$
%which is a negligible Borel set  with respect to the measure $d\rho_s(\omega)$.
It follows that for every
$\omega\in \Omega, s\in N_{\omega}^c$, $n\in\N$ and $t\in I_{\omega}$, we have
$$
1_{\{s\leq t\}} g(s,\omega)\leq \Big[g^{n,\rho}(t,\omega)+n(t-s)\Big]1_{\{s\leq t\}}.
$$
Let $\omega\in \Omega, s\in N_{\omega}^c$. If $s\in I_{\omega},$ we take $t= s$ and by letting $n$ to infinity
it follows that $g(s, \omega) \leq g^{*,\rho}(s,\omega)$. Now,
if $s\notin I_{\omega},$ then there exits a sequence $(t_p)_p\in I_{\omega}$ such that $t_p\downarrow s$. Then
$$
g(s,\omega)\leq \Big[g^{n,\rho}(t_p,\omega)+n(t_p-s)\Big].
$$
Letting $p$ to $\infty$ we have
$$
g(s,\omega)\leq g^{n,\rho}(s+,\omega) = g^{n,\rho}(s,\omega).
$$
The result follows by letting $n$ to infinity. \eop

We have the following characterization.
\begin{corollary}\label{char}
	Let $Y\in\mathcal{D}$
	\[\bigg(g(t,\omega)\leqslant Y_{t-}(\omega)\bigg) \;\; d\rho_t(\omega)P(d\omega)\text{-a.e} \quad\text{ if and only if }\quad \bigg(\forall t\in ]0,T],\;\; g^{*,\rho}(t,\omega)\leqslant Y_{t-}(\omega)\bigg)\;\; P\text{-as} \]
\end{corollary}
\bop

We use the left continuity of $Y_{t-}$ to prove that $\sup_{s\leqslant t}\big[Y_{s-}-n(t-s)\big]$ converges to $Y_{t-}$ as $n$ goes to infinity.
\eop

\end{document}